\documentclass{amsart}

\usepackage{amsthm}
\usepackage{amssymb}
\usepackage{amsmath}
\usepackage{cd}
\usepackage{epic}
\usepackage{eepic}

\newcommand     {\comment}[1]   {}
\newcommand{\mute}[2] {}
\newcommand     {\printname}[1] {}
\newcommand{\labell}[1] {\label{#1}\printname{#1}}

\DeclareMathOperator {\Orth}{Orth}

\numberwithin{equation}{section}

\newtheorem {Theorem}   {Theorem}
\newtheorem {refTheorem}   {Theorem}
\numberwithin{refTheorem}{section}
\newtheorem {Lemma}[refTheorem]{Lemma}
\newtheorem {Proposition}[refTheorem]{Proposition}

\theoremstyle{definition}

\theoremstyle{remark}
\newtheorem{Remark}[refTheorem]{Remark}
\newtheorem{Example}[refTheorem]{Example}

\newtheorem* {Claim*}{Claim}
\newtheorem* {Remark*}{Remark}

\renewcommand {\Re} {\operatorname{Re}}
\def    \valpha      {\alpha}
\def    \talpha      {\widetilde{\alpha}}

\def \bfone {{\mathbf{1}}}

\def    \cI     {{\mathcal{I}}}
\def    \calI   {{\mathcal{I}}}
\def    \cS     {{\mathcal{S}}}
\def    \calS   {{\mathcal{S}}}

\def	\calD	{{\mathcal D}}

\def	\bfC	{{\mathbf C}}
\def	\bfS	{{\mathbf S}}
\def    \bfL    {{\mathbf L}}
\def    \bfM    {{\mathbf M}}
\def    \bfT    {{\mathbf T}}

\def    \inv    {^{-1}}
\def    \ol    	{\overline}
\def    \ssminus        {{\smallsetminus}}

\def    \l<        {\left< }
\def    \r>        {\right> }
\def	\del	{\partial}
\newcommand {\deldel}[1] {\frac{\partial}{\partial #1}}

\def	\nonneg	{{\geq 0}}

\def	\shift	{{\text{shift}}}
\def    \codim    {{\operatorname{codim}}}

\def	\span	{{\operatorname{span}}}

\def	\Td	{{\operatorname{Td}}}
\def	\interior	{{\operatorname{interior}}}
\def	\relint	{{\operatorname{rel-int}}}
\def	\Vert	{{\operatorname{Vert}}}
\def	\VertDelta	{{\Vert(\Delta)}}

\def	\t	{{\mathfrak t}}

\def	\Z	{{\mathbb Z}}
\def	\R	{{\mathbb R}}
\def	\C	{{\mathbb C}}

\def	\Cinf	{C^\infty}
\def	\half	{{\frac{1}{2}}}

\def    \sss     {{\scriptstyle}}

\begin{document}

\title{Exact Euler Maclaurin formulas for simple lattice polytopes}

\author[Y.\ Karshon]{Yael Karshon}
\address{Department of Mathematics, The University of Toronto, 
Toronto, Ontario M5S 3G3, Canada}
\email{karshon@math.toronto.edu}

\author[S.\ Sternberg]{Shlomo Sternberg}
\address{Department of Mathematics, Harvard University,
Cambridge, MA 02138, USA}
\email{shlomo@math.harvard.edu}

\author[J.\ Weitsman\ \today]{Jonathan Weitsman}
\address{Department of Mathematics, University of California,
Santa Cruz, CA 95064, USA}
\email{weitsman@math.UCSC.EDU}
\thanks{This work was partially supported by
United States -- Israel Binational Science Foundation
grant number 2000352 (to Y.K. and J.W.),  by the Connaught Fund (to Y.K.),
and by National Science Foundation Grant DMS 99/71914 and 
DMS 04/05670 (to J.W.).}

\begin{abstract}
Euler Maclaurin formulas for a polytope express the sum of the values
of a function over the lattice points in the polytope in terms of
integrals of the function and its derivatives over faces of the
polytope or its expansions.   Exact Euler Maclaurin formulas
\cite{KP2,CS:bulletin,CS:EM,Gu,BV} 
apply to exponential or polynomial functions; Euler Maclaurin formulas 
with remainder \cite{KSW1,KSW2}
apply to more general smooth functions.


In this paper we review these results and present proofs of the
exact formulas obtained by these authors, using elementary methods.
We then use an algebraic formalism due to Cappell and Shaneson 
to relate the different formulas.  
\end{abstract}

\maketitle

\tableofcontents

\comment{
Cappell-Shaneson probably have formulas also for singular toric varieties.
Diaz-Robins have formulas for exponential sums over various cones.
These are probably in their papers that we're already referring to.
Look these up?

Michele said: Barvinok's recent article addresses computational aspects.
}

\section{Introduction}
\labell{sec:intro}

Let $f$ be a polynomial in one variable.
The classical Euler Maclaurin formula (see \cite[chapter XIV]{Kn})
computes the sum of the values of $f$ over the integer
points in an interval in terms of the integral of $f$ over the interval
and the values of $f$ and of its 
derivatives at the endpoints of the interval.
The formula is almost three hundred years old 
(\cite{Be}, \cite{Eu}, \cite{Mac}).
We refer the readers to the treatments by Wirtinger \cite{Wi}
and by Bourbaki \cite{Bo}.

A version of this formula that was generalized by Khovanskii and Pukhlikov 
to higher dimensions (see \cite{KP2}) involves variations of the interval.  
It reads
\begin{equation} \labell{classical EM}
f(a) + f(a+1) + \ldots + f(b-1) + f(b) =
\left. \Td \left( \deldel{h_1} \right) \Td \left( \deldel{h_2} \right)
\right|_{h_1=h_2=0}
\int_{a-h_2}^{b+h_1} f(x) dx ,
\end{equation}
where $a,b \in \Z$ and
\begin{equation} \labell{def:Td}
\Td(D) = \frac{D}{1 - e^{-D}}
       = 1 + \half D + \frac{1}{12} D^2 - \frac{1}{720} D^4 + \ldots .
\end{equation}
The right hand side of \eqref{classical EM} is well defined 
because $\int_{a-h_2}^{b+h_1} f(x) dx$, 
as a function of $h_1$ and $h_2$, is again a polynomial.
We call this an \emph{exact} formula, to distinguish it from Euler Maclaurin
formulas \emph{with remainder}, which apply to more general smooth functions.

In higher dimensions, one replaces the interval $[a,b]$ by a \emph{polytope},
that is, the convex hull of a finite set of points in a vector space, 
or, equivalently, a bounded finite intersection of closed half-spaces.  
We assume that our polytopes have nonempty interior.
A polytope in $\R^n$ is called an \emph{integral polytope},  
or a \emph{lattice polytope}, if its vertices
are in the lattice $\Z^n$.
It is called \emph{simple} if exactly $n$ edges emanate from each vertex.
For example, a two dimensional polytope (a polygon) is always simple.  
A tetrahedron and a cube are simple;
a square-based pyramid and an octahedron are not simple; 
see Figure \ref{fig:polytopes}.
\begin{figure}[ht]
\setlength{\unitlength}{0.00083333in}
\begingroup\makeatletter\ifx\SetFigFont\undefined%
\gdef\SetFigFont#1#2#3#4#5{%
  \reset@font\fontsize{#1}{#2pt}%
  \fontfamily{#3}\fontseries{#4}\fontshape{#5}%
  \selectfont}%
\fi\endgroup%
{\renewcommand{\dashlinestretch}{30}
\begin{picture}(924,939)(0,-10)
\path(12,312)(612,12)(912,312)
	(462,912)(12,312)
\path(462,912)(612,12)
\dottedline{45}(12,312)(912,312)
\end{picture}
}
\setlength{\unitlength}{0.00083333in}
\begingroup\makeatletter\ifx\SetFigFont\undefined%
\gdef\SetFigFont#1#2#3#4#5{%
  \reset@font\fontsize{#1}{#2pt}%
  \fontfamily{#3}\fontseries{#4}\fontshape{#5}%
  \selectfont}%
\fi\endgroup%
{\renewcommand{\dashlinestretch}{30}
\begin{picture}(924,939)(0,-10)
\path(12,12)(12,612)(612,612)
	(612,12)(12,12)
\path(612,12)(912,312)(912,912)(612,612)
\path(12,612)(312,912)(912,912)
\dottedline{45}(312,912)(312,312)(912,312)
\dottedline{45}(12,12)(312,312)
\end{picture}
}
\setlength{\unitlength}{0.00083333in}
\begingroup\makeatletter\ifx\SetFigFont\undefined%
\gdef\SetFigFont#1#2#3#4#5{%
  \reset@font\fontsize{#1}{#2pt}%
  \fontfamily{#3}\fontseries{#4}\fontshape{#5}%
  \selectfont}%
\fi\endgroup%
{\renewcommand{\dashlinestretch}{30}
\begin{picture}(1824,1239)(0,-10)
\dottedline{45}(12,162)(612,612)(1812,462)
\path(12,162)(1212,12)(1812,462)
\path(12,162)(912,1212)(1212,12)
\dottedline{45}(912,1212)(612,612)
\path(912,1212)(1812,462)
\end{picture}
}
\setlength{\unitlength}{0.00083333in}
\begingroup\makeatletter\ifx\SetFigFont\undefined%
\gdef\SetFigFont#1#2#3#4#5{%
  \reset@font\fontsize{#1}{#2pt}%
  \fontfamily{#3}\fontseries{#4}\fontshape{#5}%
  \selectfont}%
\fi\endgroup%
{\renewcommand{\dashlinestretch}{30}
\begin{picture}(924,939)(0,-10)
\path(12,387)(462,12)(612,312)
\path(12,387)(612,312)(912,537)
\path(462,12)(912,537)
\dottedline{45}(12,387)(312,612)(912,537)
\dottedline{45}(312,612)(462,12)
\path(12,387)(462,912)(612,312)
\path(462,912)(912,537)
\dottedline{45}(462,912)(312,612)
\end{picture}
}

\caption{}
\labell{fig:polytopes}
\end{figure}
A \emph{polytope with a non-singular fan} is a simple polytope  
in which the edges emanating from each vertex lie along vectors 
that generate the lattice $\Z^n$.
For example, in Figure \ref{fig:triangles},
the triangle on the left has a non-singular fan, and that on the right
does not (check its top vertex).

\begin{figure}[ht]
\setlength{\unitlength}{0.00083333in}
\begingroup\makeatletter\ifx\SetFigFont\undefined%
\gdef\SetFigFont#1#2#3#4#5{%
  \reset@font\fontsize{#1}{#2pt}%
  \fontfamily{#3}\fontseries{#4}\fontshape{#5}%
  \selectfont}%
\fi\endgroup%
{\renewcommand{\dashlinestretch}{30}
\begin{picture}(3376,689)(0,-10)
\put(638,637){\blacken\ellipse{40}{40}}
\put(638,637){\ellipse{40}{40}}
\put(338,637){\blacken\ellipse{40}{40}}
\put(338,637){\ellipse{40}{40}}
\put(938,337){\blacken\ellipse{40}{40}}
\put(938,337){\ellipse{40}{40}}
\put(938,37){\blacken\ellipse{40}{40}}
\put(938,37){\ellipse{40}{40}}
\put(638,337){\blacken\ellipse{40}{40}}
\put(638,337){\ellipse{40}{40}}
\put(338,337){\blacken\ellipse{40}{40}}
\put(338,337){\ellipse{40}{40}}
\put(638,37){\blacken\ellipse{40}{40}}
\put(638,37){\ellipse{40}{40}}
\put(338,37){\blacken\ellipse{40}{40}}
\put(338,37){\ellipse{40}{40}}
\put(1238,637){\blacken\ellipse{40}{40}}
\put(1238,637){\ellipse{40}{40}}
\put(1538,637){\blacken\ellipse{40}{40}}
\put(1538,637){\ellipse{40}{40}}
\put(1838,637){\blacken\ellipse{40}{40}}
\put(1838,637){\ellipse{40}{40}}
\put(2138,637){\blacken\ellipse{40}{40}}
\put(2138,637){\ellipse{40}{40}}
\put(2438,637){\blacken\ellipse{40}{40}}
\put(2438,637){\ellipse{40}{40}}
\put(2738,637){\blacken\ellipse{40}{40}}
\put(2738,637){\ellipse{40}{40}}
\put(3038,637){\blacken\ellipse{40}{40}}
\put(3038,637){\ellipse{40}{40}}
\put(3338,637){\blacken\ellipse{40}{40}}
\put(3338,637){\ellipse{40}{40}}
\put(3338,337){\blacken\ellipse{40}{40}}
\put(3338,337){\ellipse{40}{40}}
\put(3338,37){\blacken\ellipse{40}{40}}
\put(3338,37){\ellipse{40}{40}}
\put(938,637){\blacken\ellipse{40}{40}}
\put(938,637){\ellipse{40}{40}}
\put(3038,37){\blacken\ellipse{40}{40}}
\put(3038,37){\ellipse{40}{40}}
\path(1838,37)(1838,637)(3038,37)(1838,37)
\put(3038,337){\blacken\ellipse{40}{40}}
\put(3038,337){\ellipse{40}{40}}
\put(2738,337){\blacken\ellipse{40}{40}}
\put(2738,337){\ellipse{40}{40}}
\put(2438,337){\blacken\ellipse{40}{40}}
\put(2438,337){\ellipse{40}{40}}
\put(2738,37){\blacken\ellipse{40}{40}}
\put(2738,37){\ellipse{40}{40}}
\put(2438,37){\blacken\ellipse{40}{40}}
\put(2438,37){\ellipse{40}{40}}
\put(2138,337){\blacken\ellipse{40}{40}}
\put(2138,337){\ellipse{40}{40}}
\put(2138,37){\blacken\ellipse{40}{40}}
\put(2138,37){\ellipse{40}{40}}
\put(1838,337){\blacken\ellipse{40}{40}}
\put(1838,337){\ellipse{40}{40}}
\put(1838,37){\blacken\ellipse{40}{40}}
\put(1838,37){\ellipse{40}{40}}
\put(1538,337){\blacken\ellipse{40}{40}}
\put(1538,337){\ellipse{40}{40}}
\put(1538,37){\blacken\ellipse{40}{40}}
\put(1538,37){\ellipse{40}{40}}
\put(1238,37){\blacken\ellipse{40}{40}}
\put(1238,37){\ellipse{40}{40}}
\put(1238,337){\blacken\ellipse{40}{40}}
\put(1238,337){\ellipse{40}{40}}
\put(38,637){\blacken\ellipse{40}{40}}
\put(38,637){\ellipse{40}{40}}
\put(38,337){\blacken\ellipse{40}{40}}
\put(38,337){\ellipse{40}{40}}
\put(38,37){\blacken\ellipse{40}{40}}
\put(38,37){\ellipse{40}{40}}
\path(338,37)(338,637)(938,37)(338,37)
\end{picture}
}
\caption{}
\labell{fig:triangles}
\end{figure}
We refer the reader to \cite{B,Gr,KW,Z} for general background on convex 
polytopes.

\begin{Remark}
A \emph{fan} in $\R^n$ is a set of convex polyhedral cones emanating
from the origin, such that the intersection of any two cones in the set 
is a common face, and such that the face of any cone in the set is itself
a cone in the set.  The fan of a polytope $\Delta \subset \R^n$ consists
of a set of cones associated to the faces of $\Delta$;
for each face we take the cone
generated by the inward normals to the facets that meet at that face.

A convex polyhedral cone is \emph{non-singular} if it can be generated
by a set of vectors in $\Z^n$ which are part of a $\Z$-basis of $\Z^n$.
A fan is non-singular if each cone in the fan is non-singular.
The name ``non-singular" comes from the theory of toric varieties;
non-singular fans correspond to non-singular toric varieties.

The terminology in the literature is inconsistent.  A non-singular cone
is also called ``smooth cone" and ``unimodular cone". 
Polytopes with non-singular fans are also called ``non-singular polytopes",
``smooth polytopes", or ``Delzant polytopes".  
In our previous papers \cite{KSW1,KSW2} we used the terms ``regular orthant" 
and ``regular polytope" (not to be confused with the more common usage 
of this term as ``platonic solid").
Other terms that have been suggested to us are ``unimodular polytope"
or ``torsion-free polytope".  
\end{Remark}

Khovanskii and Pukhlikov \cite{KP1,KP2}, following
Khovanskii \cite{Kh1,Kh2} (see also Kantor and Khovanskii \cite{KK,KK2}),
gave a formula
for the sum of the values of a quasi-polynomial (polynomial times
exponential) function on the lattice points in a lattice polytope 
with non-singular fan.
This formula was further generalized to \emph{simple} polytopes 
by Cappell and Shaneson \cite{CS:bulletin,CS:EM,CS:private,S},
and subsequently by Guillemin \cite{Gu} and by Brion-Vergne 
\cite{BV,BV:partition}.
Cappell-Shaneson and Brion-Vergne \cite{BV:partition} also work
with polytopes that are not simple; in this paper we will restrict
ourselves to simple polytopes. 
Also see the explicit formulas in \cite{SV}
and the survey \cite{V:survey}.
When applied to the constant function $f \equiv 1$, these formulas
compute the number of lattice points in a simple lattice polytope $\Delta$
in terms of the volumes of ``expansions" of $\Delta$.
A sample of the literature on the problem of counting lattice
points in convex polytopes is \cite{pick,Md,V,KK,Mo,Po,DR,BDR,Ha,dCP};
see the survey \cite{BP} and references therein.

\begin{Remark} \labell{alg-geom}
Khovanskii's motivation came from algebraic geometry: a lattice polytope 
$\Delta$ with non-singular fan determines a toric variety $M_\Delta$
and a holomorphic line bundle $\bfL_\Delta \to M_\Delta$.
The quantization $Q(M_\Delta)$ is interpreted as the
space of holomorphic sections of $L_\Delta$ and is computed
by the Hirzebruch Riemann Roch formula.
The lattice points in $\Delta$ correspond to basis elements
of $Q(M_\Delta)$.
A \emph{simple} polytope $\Delta$ still determines a toric variety, 
which now may have orbifold singularities.  
Cappell and Shaneson derived their formula from their
theory of characteristic classes of singular algebraic varieties.
Guillemin derived his formula by applying the Kawasaki-Riemann-Roch formula
to symplectic toric orbifolds.
Brion and Vergne's proof uses Fourier analysis and is closer to Khovanskii
and Pukhlikov's original proof.
\end{Remark}

In this paper we present an elementary proof of the exact 
Euler-Maclaurin formulas that follows the lines of the original 
Khovanskii-Pukhlikov proof, through a decomposition of the polytope
into an alternating sum of simple convex polyhedral cones.
We then use an algebraic formalism due to Cappell and Shaneson
to explain the equivalence of the different formulas.

The proof of the exact Euler Maclaurin formula for a simple convex
polyhedral cone involves the following ingredients: 
the summation of a geometric series, the change of variable
formula for integration, and Frobenius's theorem that the average value
of a non-trivial character of a finite group is zero.
(See section \ref{sec:first formula}.)
The ``polar decomposition" of the polytope into simple convex polyhedral cones
was proved in papers of Varchenko and Lawrence \cite{lawrence,V}. 
We present a short direct proof of it.  (See section \ref{sec:decompose}.)

From a simple polytope $\Delta$ with $d$ faces one gets 
a expanded polytope $\Delta(h)$, for $h=(h_1,\ldots,h_d)$,
by parallel translating the hyperplanes containing the facets,
see equation \eqref{polytope} below.
The integrals of a function $f$ on $\Delta(h)$ and on its faces
are functions of the $d$ variables $h_1,\ldots,h_d$.
The formulas of Khovanskii-Pukhlikov, Guillemin, and Brion-Vergne
involve an application of infinite order differential operators
to these functions.
The Cappell-Shaneson formula does not involve expansions of the polytope.  
It is stated through a formalism that we call the 
Cappell-Shaneson algebra.  Their abstract formula translates 
to several different concrete formulas; each of these involves 
applying differential operators to the function and integrating over 
faces of the polytope.  
The relations in the Cappell-Shaneson algebra allow one to pass 
between the different concrete formulas.  In Section \ref{sec:stokes} 
we show how to incorporate expansions in $h$ into the Cappell-Shaneson 
formalism.
In section \ref{sec:CS formula} we use a generalization of the 
``polar decomposition",
which applies to polytopes with some facets removed, to prove that
the Cappell-Shaneson formula is equivalent to the Khovanskii-Pukhlikov 
formula in the case of polytopes with non-singular fans.

\smallskip\noindent\textbf{Acknowledgement}
We would like to thank Itai Bar-Natan, S. Cappell, V. Guillemin, G. Kalai, 
G. Kuperberg, and J. Shaneson for helpful discussions.
We are particularly grateful to A. Khovanskii and to 
M. Vergne for explaining to us their proofs.

\section{Euler Maclaurin formulas in one dimension}
\labell{sec:Todd}

In this section we present Euler-MacLaurin formulas 
for a ray and for an interval, in order to illustrate arguments
that generalize to higher dimensions.


\smallskip\noindent\textbf{The ODE for the exponential function.}\ 
Let $D = \displaystyle{\deldel{h}}$.  Since
$$D^{k}e^{\xi h} = \xi^{k}e^{\xi h},$$
for any formal power series $F$ in one variable we have
\begin{equation} \labell{generalrule}
F(D)e^{\xi h}= F(\xi)e^{\xi h}
\end{equation}
in the ring of power series in two variables.
It follows that, for any non-negative integer $N$,
\begin{equation} \labell{elephant}
 F(D) \frac{ (\xi h)^N }{N!} = \left( F(\xi) e^{\xi h} \right)^{\l< N \r> }
\end{equation}
where the superscript $\l< N \r> $ denotes the $N$th term in the Taylor
expansion in $\xi$.  Under suitable convergence conditions,
\eqref{generalrule} is an  equality of functions.  See \cite{Bo}. 

\smallskip\noindent\textbf{Euler Maclaurin formula for a ray.}\ 
To conform with the topological literature, let us define
the \textbf{Todd function} by
\begin{equation} \labell{Toddfunction}
\Td(S) := \frac{S}{1-e^{-S}}
\end{equation}
and the corresponding \textbf{Todd operator} in one variable by
$$
\Td\left( D \right).
$$

Our general rule (\ref{generalrule}) gives
$$\Td\left( D \right)e^{\xi h}=\Td(\xi)e^{\xi h}$$
in the ring of formal power series.
If $|\xi| < 2\pi$, so that the power series for $\Td(\xi) $ converges,
we can regard this last equation as an equality of functions. 
Namely, the left hand side is the limit of the functions
obtained by applying the partial sums of the infinite series 
$\Td( D )$ to the exponential function.
If $\xi\neq 0$, we can re-write this as
\begin{equation} \labell{Td eq}
\Td\left( D \right) \frac{e^{\xi h}}{\xi}
  = \frac {e^{\xi h}} {1-e^{-\xi}}.
\end{equation}

If $\xi > 0$, the geometric series expansion
$$ 1+e^{-\xi}+e^{-2\xi}+e^{-3\xi}+\cdots = \frac1{1-e^{-\xi}} $$
converges, as does the integral
$$\int_{-\infty}^h e^{\xi x}dx = \frac{e^{\xi h}}{\xi},$$
so \eqref{Td eq} gives
\begin{equation} \labell{EM ray}
  \left.\Td\left(\deldel{h}\right)
  \int_{-\infty}^h e^{\xi x} dx \right|_{h=0}
  = \sum_{n=-\infty}^0 e^{\xi n}.
\end{equation}
This is the Euler Maclaurin formula for the ray $(-\infty,0]$,
with the function $f(x) = e^{\xi x}$.

\smallskip\noindent\textbf{Polar decomposition of an interval.}\ 
In the one dimensional case, the ``polar decomposition" becomes
the relation 
\begin{equation}  \labell{Lawrence for interval}
 \mathbf{1}_{I}(x)
   = \mathbf{1}_{\bfC_b}(x) - \mathbf{1}_{\bfC_a^\sharp}(x)  
\end{equation}
between the characteristic functions of the interval
$I = [a,b]$, the ray $\bfC_b = (-\infty,b]$, and 
the ray $\bfC_a^\sharp = (-\infty,a)$
(which is obtained from the ray $C_a = [a,\infty)$
by flipping its direction and removing its vertex).

\smallskip\noindent\textbf{Euler Maclaurin on finite intervals.}\ 
Let $I=[a,b]$ be a closed interval with integer endpoints.
For $h = (h_1,h_2)$, consider the expanded interval
$$ I(h) := [a-h_2,b+h_1]. $$
Summation and integration of the function
$$ f(x) = e^{\xi x} $$
gives
\begin{equation} \labell{five}
  \calI(h,\xi) := \int_{I(h)} e^{\xi x} dx 
  = \frac{e^{\xi(b+h_1)}}{\xi} - \frac{e^{\xi(a-h_2)}}{\xi}
\end{equation}
for all $\xi$ such that $\xi \neq 0$, and
\begin{equation} \labell{ten}
  \calS(\xi) := \sum_{x \in I \cap \Z} e^{\xi x} 
  = \frac{e^{\xi b}}{1 - e^{-\xi}} + \frac{e^{\xi a}}{1 - e^{\xi}}
\end{equation}
for all $\xi \in \C$ such that $e^{\xi} \neq 1$.
An indirect proof of \eqref{five} and \eqref{ten},
which generalizes to higher dimensions, uses 
the ``polar decomposition" \eqref{Lawrence for interval}:  
if $\Re \xi > 0$, then
$$ \sum_{k=-\infty}^b e^{\xi x} dx = \frac{e^{\xi b}}{1 - e^{-\xi}} \ 
\quad \text{ and } \quad
   \sum_{k=-\infty}^{a-1} e^{\xi x} dx = \frac{e^{\xi (a-1)}}{1 - e^{-\xi}}
   = - \frac{e^{\xi a}}{1 - e^\xi} \ . $$
Since $\calS(\xi)$ is the difference of these two infinite sums,
\eqref{ten} holds whenever $\Re \xi > 0$.
Because the set $\{ \xi \in \C \ | \ e^\xi \neq 1 \}$ is connected,
by analytic continuation \eqref{ten} holds for \emph{all} $\xi$ in this set.
A similar argument shows that \eqref{five} holds for all $\xi$ in the set 
$\{ \xi \in \C \ | \ \xi \neq 0 \}$.

\medskip

At this point one can proceed in several ways.

\medskip

\smallskip\noindent\textbf{Formal approach.}\ 
One can deduce an Euler Maclaurin formula for polynomial functions
directly from \eqref{five} and \eqref{ten}.
This is the one dimensional case of the approach of Brion-Vergne.
From \eqref{five} we get
\begin{equation} \labell{six}
\xi \calI(h,\xi) = e^{\xi(b+h_1)} - e^{\xi(a-h_2)}
\end{equation}
for all $\xi \neq 0$, and, by continuity, also for $\xi = 0$.
From \eqref{ten} we get
\begin{equation} \labell{seven}
\xi \calS(\xi) = \Td(\xi) e^{\xi b} - \Td(-\xi) e^{\xi a}
\end{equation}
for all $\xi$ such that $e^\xi \neq 1$ and, by continuity,
also for $\xi = 0$.  Comparing the Taylor coefficients with respect to $\xi$
on the left and right hand sides of \eqref{six} and of \eqref{seven},
we get
\begin{equation} \labell{Hee Hee}
   \xi \int_{I(h)} \frac{(\xi x)^N}{N!} dx
   = \frac{(\xi(b+h_1))^{N+1}}{(N+1)!} 
   - \frac{(\xi(a-h_2))^{N+1}}{(N+1)!}
\end{equation}
and
\begin{equation} \labell{Hoo Hoo}
\xi \sum_{x \in I \cap \Z} \frac{(\xi x)^N}{N!}
 = (e^{\xi b} \Td(\xi) - e^{\xi a} \Td(-\xi)  )^{\langle N+1 \rangle },
\end{equation}
where the superscript $\langle N+1 \rangle $ denotes the summand
that is homogeneous of degree \ $N+1$ \ in $\xi$.  
Since $\displaystyle{\Td(\deldel{h_i})} = $ 1 + a multiple of 
$\displaystyle{\deldel{h_i}}$,
\begin{multline*}
   \left. \Td\left(\deldel{h_1}\right)
          \Td\left(\deldel{h_2}\right) \right|_{h=0}
   \xi \int_{I(h)} \frac{(\xi x)^N}{N!} dx \\
   =  \left. \Td\left(\deldel{h_1}\right)\right|_{h_1 = 0}
      \frac{(\xi(b+h_1))^{N+1}}{(N+1)!}
   -  \left. \Td\left(\deldel{h_2}\right)\right|_{h_2 = 0}
      \frac{(\xi(a-h_2))^{N+1}}{(N+1)!} \qquad \text{by \eqref{Hee Hee}}, \\
 = \left( \Td(\xi) e^{\xi b} - \Td(-\xi) e^{\xi a} \right)
                                                ^{ \langle N+1 \rangle} 
\qquad \text{by \eqref{elephant}}, \\
   = \xi \sum_{x \in I \cap \Z} \frac{(\xi x)^N}{N!}  
\qquad \text{by \eqref{Hoo Hoo}}.
\end{multline*}
This gives the Euler-Maclaurin formula
\begin{equation} \labell{ultimate I}
 \left. \Td(\deldel{h_1}) \Td(\deldel{h_2}) \right|_{h=0} 
 \int_{I(h)} f = \sum_{x \in I \cap \Z} f
\end{equation}
for the function $\displaystyle{f(x) = \frac{\xi^{N+1} x^N}{N!}}$,
whenever $\xi \neq 0$.
Because multiplication by $N!$ and division by the non-zero constant 
$\xi^{N+1}$ commutes with summation, with integration, and with the
infinite order differential operator $\Td(\deldel{h_1}) \Td(\deldel{h_2})$,
we deduce the Euler-Maclaurin formula \eqref{ultimate I}
for the monomials $f(x) = x^N$, and hence for all polynomials.

\medskip

\smallskip\noindent\textbf{Approach through Euler Maclaurin 
for exponentials.}\ 
In the original approach of Khovanskii-Pukhlikov,
one deduces an Euler Maclaurin formula for polynomials, and, 
more generally, for (quasi)-polynomials, from a formula for exponentials.  
(A quasi-polynomial is a sum of products of exponentials by polynomials.)
In the one dimensional case,
the Euler-Maclaurin formula for exponentials asserts that
\begin{equation} \labell{EM exp}
\left. \Td\left(\deldel{h_1}\right) \Td\left(\deldel{h_2}\right)
\right|_{h=0} \int_{I(h)} e^{\xi x} = \sum_{I \cap \Z} e^{\xi x},
\end{equation}
or, equivalently, that
$$ \left. \Td\left(\deldel{h_1}\right) \Td\left(\deldel{h_2}\right) 
   \right|_{h=0}
   \calI(h,\xi) = \calS(\xi) .$$
This formula is true for all $\xi$ such that $|\xi| < 2\pi$.
For $\xi$ in the punctured disk
\begin{equation} \labell{the set}
 \{ \xi \in \C \ | \ \xi \neq 0 \ , \ |\xi| < 2\pi \},
\end{equation}
the formula follows immediately from \eqref{Lawrence for interval},
\eqref{five}, \eqref{ten}, and from the facts that
\begin{equation} \labell{face-a}
 \left. \Td\left(\deldel{h_1}\right) \right|_{h=0} \frac{e^{\xi(b+h_1)}}{\xi}
 = \Td(\xi) \frac{e^{\xi b}}{\xi} = \frac{e^{\xi b}}{1 - e^{-\xi}}
\end{equation}
and
\begin{equation} \labell{face-b}
 \left. \Td\left(\deldel{h_2}\right) \right|_{h=0} \frac{e^{\xi(a-h_2)}}{\xi}
 = \Td(-\xi) \frac{e^{\xi a}}{\xi} = - \frac{e^{\xi a}}{1 - e^{\xi}}.
\end{equation}
(If $\Re \xi \neq 0$ then \eqref{face-a} is an Euler-Maclaurin formula 
for the ray $(-\infty,b]$ or $[b,\infty)$, and similarly for \eqref{face-b}.
However, \eqref{face-a} and \eqref{face-b} hold for \emph{all} $\xi$
in the set \eqref{the set}.)

In \eqref{face-a} and \eqref{face-b}, the left hand sides converge 
to the right hand sides uniformly in $\xi$ on compact subsets of the 
punctured disk \eqref{the set}.
This is because the Taylor series of $\Td(\cdot)$ converges uniformly
on compact subsets of the disk $\{ |\xi| < 2\pi \}$, and the functions
$\frac{e^{\xi b}}{\xi}$ and $\frac{e^{\xi a}}{\xi}$ are bounded
on compact subsets of the punctured disk \eqref{the set}.
It follows that in \eqref{EM exp} the left hand side converges 
to the right hand side uniformly in $\xi$ on compact subsets 
of \eqref{the set}.  But the right hand side and the partial sums
of the left hand side of \eqref{EM exp} are analytic in $\xi$
for \emph{all} $|\xi| < 2\pi$.  It follows from the Cauchy integral formula
that the left hand side of \eqref{EM exp} converges to the right hand side,
uniformly on compact subsets, on \emph{all} of
$\{ \xi \ | \ |\xi| < 2\pi \}$.
 
It further follows that the infinite sum on the
left hand side of \eqref{EM exp} can be differentiated with respect to $\xi$
term by term.  Hence, the infinite order differential operator on the
left hand side of \eqref{EM exp} commutes with differentiation with respect
to $\xi$.  
Since 
$$ \frac{\del^k}{\del \xi^k} \int_{I(h)} e^{\xi x} dx
   = \int_{I(h)} x^k e^{\xi x} dx 
\qquad \text{ and } \qquad
\frac{\del^k}{\del \xi^k} \sum_{I \cap \Z} e^{\xi x}
\; = \; \sum_{I \cap \Z} x^k e^{\xi x} ,
$$
we get the Euler-Maclaurin formula \eqref{ultimate I}
for the function $f(x) = x^k e^{\xi x}$ by differentiating
the left and right hand sides of \eqref{EM exp} $k$ times
with respect to $\xi$.  

\medskip

\smallskip\noindent\textbf{Approach through Euler Maclaurin for rays.}\ 
Yet another approach is to deduce the Euler Maclaurin formula 
for an interval directly from the Euler Maclaurin formula for a ray.
See Appendix \ref{sec:polynomials}.

\section{Polar decomposition of a simple polytope}
\labell{sec:decompose}

\begin{figure}[ht]
\begin{center}
\setlength{\unitlength}{0.0004in}
\begingroup\makeatletter\ifx\SetFigFont\undefined%
\gdef\SetFigFont#1#2#3#4#5{%
  \reset@font\fontsize{#1}{#2pt}%
  \fontfamily{#3}\fontseries{#4}\fontshape{#5}%
  \selectfont}%
\fi\endgroup%
{\renewcommand{\dashlinestretch}{30}
\begin{picture}(11649,3039)(0,-10)
\texture{0 0 0 888888 88000000 0 0 80808
        8000000 0 0 888888 88000000 0 0 80808
        8000000 0 0 888888 88000000 0 0 80808
        8000000 0 0 888888 88000000 0 0 80808 }
%
%
\shade\path(12,3012)(12,1812)(1212,1812)(12,3012)
\put(12,3012){\blacken\ellipse{60}{60}}
\put(12,2712){\blacken\ellipse{60}{60}}
\put(312,2712){\blacken\ellipse{60}{60}}
\put(12,2412){\blacken\ellipse{60}{60}}
\put(312,2412){\blacken\ellipse{60}{60}}
\put(612,2412){\blacken\ellipse{60}{60}}
\put(12,2112){\blacken\ellipse{60}{60}}
\put(312,2112){\blacken\ellipse{60}{60}}
\put(612,2112){\blacken\ellipse{60}{60}}
\put(912,2112){\blacken\ellipse{60}{60}}
\put(12,1812){\blacken\ellipse{60}{60}}
\put(312,1812){\blacken\ellipse{60}{60}}
\put(612,1812){\blacken\ellipse{60}{60}}
\put(912,1812){\blacken\ellipse{60}{60}}
\put(1212,1812){\blacken\ellipse{60}{60}}
%
%
\put(1212,2312){$=$}
\shade\path(1812,3012)(1812,12)(2487,237)
	(3387,162)(3387,162)(3987,312)
	(4512,312)(1812,3012)
\dottedline{45}(1812,1812)(4212,1812)
\put(1812,3012){\blacken\ellipse{60}{60}}
\put(1812,2712){\blacken\ellipse{60}{60}}
\put(2112,2712){\blacken\ellipse{60}{60}}
\put(1812,2412){\blacken\ellipse{60}{60}}
\put(2112,2412){\blacken\ellipse{60}{60}}
\put(2412,2412){\blacken\ellipse{60}{60}}
\put(1812,2112){\blacken\ellipse{60}{60}}
\put(2112,2112){\blacken\ellipse{60}{60}}
\put(2412,2112){\blacken\ellipse{60}{60}}
\put(2712,2112){\blacken\ellipse{60}{60}}
\put(1812,1812){\blacken\ellipse{60}{60}}
\put(2112,1812){\blacken\ellipse{60}{60}}
\put(2412,1812){\blacken\ellipse{60}{60}}
\put(2712,1812){\blacken\ellipse{60}{60}}
\put(3012,1812){\blacken\ellipse{60}{60}}
\put(1812,1512){\blacken\ellipse{60}{60}}
\put(2112,1512){\blacken\ellipse{60}{60}}
\put(2412,1512){\blacken\ellipse{60}{60}}
\put(2712,1512){\blacken\ellipse{60}{60}}
\put(3012,1512){\blacken\ellipse{60}{60}}
\put(3312,1512){\blacken\ellipse{60}{60}}
\put(1812,1212){\blacken\ellipse{60}{60}}
\put(2112,1212){\blacken\ellipse{60}{60}}
\put(2412,1212){\blacken\ellipse{60}{60}}
\put(2712,1212){\blacken\ellipse{60}{60}}
\put(3012,1212){\blacken\ellipse{60}{60}}
\put(3312,1212){\blacken\ellipse{60}{60}}
\put(3612,1212){\blacken\ellipse{60}{60}}
\put(1812,912){\blacken\ellipse{60}{60}}
\put(2112,912){\blacken\ellipse{60}{60}}
\put(2412,912){\blacken\ellipse{60}{60}}
\put(2712,912){\blacken\ellipse{60}{60}}
\put(3012,912){\blacken\ellipse{60}{60}}
\put(3312,912){\blacken\ellipse{60}{60}}
\put(3612,912){\blacken\ellipse{60}{60}}
\put(3912,912){\blacken\ellipse{60}{60}}
\put(1812,612){\blacken\ellipse{60}{60}}
\put(2112,612){\blacken\ellipse{60}{60}}
\put(2412,612){\blacken\ellipse{60}{60}}
\put(2712,612){\blacken\ellipse{60}{60}}
\put(3012,612){\blacken\ellipse{60}{60}}
\put(3312,612){\blacken\ellipse{60}{60}}
\put(3612,612){\blacken\ellipse{60}{60}}
\put(3912,612){\blacken\ellipse{60}{60}}
\put(4212,612){\blacken\ellipse{60}{60}}
\put(1812,312){\blacken\ellipse{60}{60}}
\put(2112,312){\blacken\ellipse{60}{60}}
\put(2412,312){\blacken\ellipse{60}{60}}
\put(2712,312){\blacken\ellipse{60}{60}}
\put(3012,312){\blacken\ellipse{60}{60}}
\put(3312,312){\blacken\ellipse{60}{60}}
\put(3612,312){\blacken\ellipse{60}{60}}
\put(1812,12){\blacken\ellipse{60}{60}}
%
%
\put(4312,2312){$-$}
\shade\path(5112,1812)(5112,162)(5487,312)
	(6762,312)(7137,612)(7662,1437)
	(7812,1812)(5112,1812)
\dottedline{45}(5112,1812)(5112,3012)(6312,1812)(7262,862)
\put(5112,1512){\blacken\ellipse{60}{60}}
\put(5412,1512){\blacken\ellipse{60}{60}}
\put(5712,1512){\blacken\ellipse{60}{60}}
\put(6012,1512){\blacken\ellipse{60}{60}}
\put(6312,1512){\blacken\ellipse{60}{60}}
\put(6612,1512){\blacken\ellipse{60}{60}}
\put(6912,1512){\blacken\ellipse{60}{60}}
\put(7212,1512){\blacken\ellipse{60}{60}}
\put(7512,1512){\blacken\ellipse{60}{60}}
\put(5112,1212){\blacken\ellipse{60}{60}}
\put(5412,1212){\blacken\ellipse{60}{60}}
\put(5712,1212){\blacken\ellipse{60}{60}}
\put(6012,1212){\blacken\ellipse{60}{60}}
\put(6312,1212){\blacken\ellipse{60}{60}}
\put(6612,1212){\blacken\ellipse{60}{60}}
\put(6912,1212){\blacken\ellipse{60}{60}}
\put(7212,1212){\blacken\ellipse{60}{60}}
\put(5112,912){\blacken\ellipse{60}{60}}
\put(5412,912){\blacken\ellipse{60}{60}}
\put(5712,912){\blacken\ellipse{60}{60}}
\put(6012,912){\blacken\ellipse{60}{60}}
\put(6312,912){\blacken\ellipse{60}{60}}
\put(6612,912){\blacken\ellipse{60}{60}}
\put(6912,912){\blacken\ellipse{60}{60}}
\put(7212,912){\blacken\ellipse{60}{60}}
\put(5112,612){\blacken\ellipse{60}{60}}
\put(5412,612){\blacken\ellipse{60}{60}}
\put(5712,612){\blacken\ellipse{60}{60}}
\put(6012,612){\blacken\ellipse{60}{60}}
\put(6312,612){\blacken\ellipse{60}{60}}
\put(6612,612){\blacken\ellipse{60}{60}}
\put(6912,612){\blacken\ellipse{60}{60}}
\put(5112,312){\blacken\ellipse{60}{60}}
\put(5412,312){\blacken\ellipse{60}{60}}
%
%
\put(7812,2312){$+$}
\shade\path(9612,1812)(11412,12)(11562,912)
	(11462,1362)(11437,1812)(9612,1812)
\dottedline{45}(8412,3012)(8412,1812)(9612,1812)(8412,3012)
\dottedline{45}(8412,1812)(8412,112)
\put(10212,1512){\blacken\ellipse{60}{60}}
\put(10512,1512){\blacken\ellipse{60}{60}}
\put(10812,1512){\blacken\ellipse{60}{60}}
\put(11112,1512){\blacken\ellipse{60}{60}}
\put(11412,1512){\blacken\ellipse{60}{60}}
\put(10512,1212){\blacken\ellipse{60}{60}}
\put(10812,1212){\blacken\ellipse{60}{60}}
\put(11112,1212){\blacken\ellipse{60}{60}}
\put(11412,1212){\blacken\ellipse{60}{60}}
\put(10812,912){\blacken\ellipse{60}{60}}
\put(11112,912){\blacken\ellipse{60}{60}}
\put(11412,912){\blacken\ellipse{60}{60}}
\put(11112,612){\blacken\ellipse{60}{60}}
\put(11412,612){\blacken\ellipse{60}{60}}
\put(11412,312){\blacken\ellipse{60}{60}}
\end{picture}
}
\end{center}
\caption{The polar decomposition theorem}
\label{fig:Lawrence}
\end{figure}
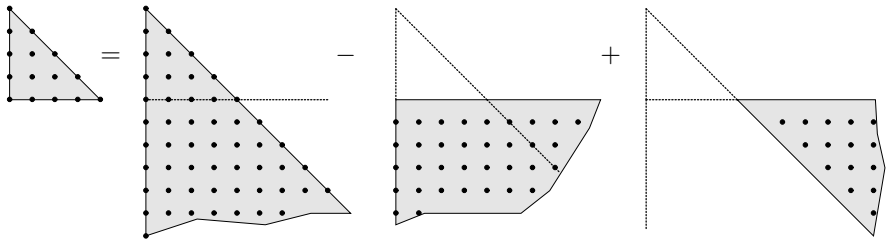

In this section we describe a decomposition 
of a simple polytope into simple convex polyhedral cones.  These cones
have apexes at the vertices of the polytope.  Each is generated 
by flipping some of the edge vectors according to a choice of 
``polarization", so that they all point roughly in the same direction,
and removing corresponding facets.  
For an illustration of this decomposition in the case of a triangle,
see Figure \ref{fig:Lawrence}.

In this section we present a short direct proof of the polar
decomposition of a simple polytope, similar to the one that we gave 
in \cite{KSW2}.  In Section \ref{sec:CS formula} 
(see \eqref{Lawrence without facets}) we give a variant
of this decomposition that applies to a polytope with some facets removed. 

Let $\Delta$ be a polytope in an $n$ dimensional vector space $V$
and $F$ a face of $\Delta$.
The \emph{tangent cone} to $\Delta$ at $F$ is
$$ \bfC_F =
   \{ y + r (x-y) \ | \ r \geq 0 \, , \, y \in F \, , \, x \in \Delta \} .$$
(Warning: other authors define the tangent cone as 
$\{ r (x-y) \ | \ r \geq 0 \, , \, y \in F \, , \, x \in \Delta \}$.)

\begin{figure}[ht]
\setlength{\unitlength}{0.0006in}
\begingroup\makeatletter\ifx\SetFigFont\undefined%
\gdef\SetFigFont#1#2#3#4#5{%
  \reset@font\fontsize{#1}{#2pt}%
  \fontfamily{#3}\fontseries{#4}\fontshape{#5}%
  \selectfont}%
\fi\endgroup%
{\renewcommand{\dashlinestretch}{30}
\begin{picture}(2529,1826)(0,-10)
\thinlines
\texture{0 0 0 888888 88000000 0 0 80808 
	8000000 0 0 888888 88000000 0 0 80808 
	8000000 0 0 888888 88000000 0 0 80808 
	8000000 0 0 888888 88000000 0 0 80808 }
\shade\path(105,1799)(2505,599)
        (2506,596)
	(2507,592)(2508,587)(2510,579)
	(2512,570)(2514,559)(2516,546)
	(2517,532)(2517,517)(2516,501)
	(2513,485)(2509,467)(2503,450)
	(2493,431)(2481,411)(2465,391)
	(2445,369)(2420,347)(2390,323)
	(2355,299)(2322,279)(2289,260)
	(2256,243)(2226,228)(2199,215)
	(2176,203)(2156,194)(2140,186)
	(2127,179)(2116,174)(2108,169)
	(2100,165)(2092,161)(2085,158)
	(2076,154)(2065,151)(2051,146)
	(2033,141)(2011,135)(1983,128)
	(1948,121)(1908,112)(1860,103)
	(1805,93)(1745,83)(1680,74)
	(1626,67)(1572,62)(1519,56)
	(1470,52)(1424,48)(1382,44)
	(1344,41)(1310,38)(1280,35)
	(1253,32)(1230,29)(1209,27)
	(1190,24)(1174,22)(1158,20)
	(1143,18)(1127,16)(1112,14)
	(1095,13)(1077,12)(1056,12)
	(1034,12)(1008,13)(979,14)
	(946,17)(910,21)(870,26)
	(826,32)(779,40)(730,49)
	(680,61)(630,74)(575,91)
	(525,109)(479,127)(440,145)
	(406,161)(377,175)(354,188)
	(335,200)(320,209)(309,218)
	(301,225)(295,231)(290,237)
	(286,243)(283,249)(279,255)
	(275,263)(270,272)(263,284)
	(254,297)(243,314)(229,333)
	(212,356)(194,383)(173,413)
	(150,447)(127,485)(105,524)
	(83,568)(65,610)(49,650)
	(37,687)(27,719)(20,747)
	(16,770)(13,790)(12,806)
	(12,819)(13,830)(15,840)
	(18,849)(20,858)(23,868)
	(26,880)(28,893)(30,910)
	(31,931)(32,956)(33,985)
	(33,1020)(32,1059)(31,1103)
	(30,1150)(30,1199)(30,1252)
	(32,1302)(34,1349)(37,1393)
	(40,1434)(44,1472)(48,1508)
	(52,1541)(57,1572)(62,1602)
	(67,1630)(72,1656)(77,1681)
	(82,1704)(87,1725)(91,1743)
	(95,1759)(98,1773)(101,1783)
	(103,1790)(104,1795)(105,1798)(105,1799)
\thicklines
\path(105,1799)(2505,599)
\path(705,1499)(1905,899)(705,899)(705,1499)
\end{picture}
}
\setlength{\unitlength}{0.0006in}
\begingroup\makeatletter\ifx\SetFigFont\undefined%
\gdef\SetFigFont#1#2#3#4#5{%
  \reset@font\fontsize{#1}{#2pt}%
  \fontfamily{#3}\fontseries{#4}\fontshape{#5}%
  \selectfont}%
\fi\endgroup%
{\renewcommand{\dashlinestretch}{30}
\begin{picture}(2498,1735)(0,-10)
\texture{0 0 0 888888 88000000 0 0 80808
        8000000 0 0 888888 88000000 0 0 80808
        8000000 0 0 888888 88000000 0 0 80808
        8000000 0 0 888888 88000000 0 0 80808 }
\thinlines
\shade\path(33,487)(33,1687)(2433,487)
        (2434,486)(2435,484)
	(2437,481)(2440,476)(2444,470)
	(2449,461)(2454,451)(2460,440)
	(2466,427)(2472,413)(2478,398)
	(2482,383)(2485,367)(2486,351)
	(2486,334)(2483,317)(2477,300)
	(2468,282)(2455,264)(2438,245)
	(2417,226)(2390,206)(2358,187)
	(2325,170)(2290,155)(2256,141)
	(2224,129)(2195,119)(2170,111)
	(2148,105)(2129,99)(2114,95)
	(2101,92)(2090,90)(2080,88)
	(2070,87)(2061,86)(2049,84)
	(2036,83)(2019,81)(1998,79)
	(1972,75)(1941,72)(1902,67)
	(1857,62)(1804,56)(1745,50)
	(1679,43)(1608,37)(1549,33)
	(1490,29)(1434,26)(1380,23)
	(1329,20)(1283,18)(1241,17)
	(1203,16)(1169,15)(1139,14)
	(1112,13)(1088,13)(1066,12)
	(1045,12)(1026,12)(1008,12)
	(990,12)(971,12)(952,12)
	(931,13)(908,13)(883,14)
	(855,15)(824,16)(790,17)
	(752,18)(712,20)(669,23)
	(623,26)(576,29)(529,33)
	(483,37)(422,44)(370,52)
	(268,98)(267,104)(262,112)
	(253,122)(240,134)(223,150)
	(204,167)(183,187)(159,213)
	(138,238)(120,264)(105,288)
	(92,311)(81,334)(71,356)
	(63,377)(56,398)(50,417)
	(45,435)(41,451)(38,464)
	(36,474)(34,481)(33,485)(33,487)
\thicklines
\path(33,1687)(1233,1087)(33,1087)(33,1687)
\path(33,487)(33,1687)(2433,487)
\end{picture}
}
\setlength{\unitlength}{0.0006in}
\begingroup\makeatletter\ifx\SetFigFont\undefined%
\gdef\SetFigFont#1#2#3#4#5{%
  \reset@font\fontsize{#1}{#2pt}%
  \fontfamily{#3}\fontseries{#4}\fontshape{#5}%
  \selectfont}%
\fi\endgroup%
{\renewcommand{\dashlinestretch}{30}
\begin{picture}(2441,1677)(0,-10)
\thinlines
\texture{0 0 0 888888 88000000 0 0 80808 
	8000000 0 0 888888 88000000 0 0 80808 
	8000000 0 0 888888 88000000 0 0 80808 
	8000000 0 0 888888 88000000 0 0 80808 }
\shade\path(239,1587)(278,1599)(318,1609)
	(357,1617)(395,1624)(431,1630)
	(464,1635)(494,1639)(522,1642)
	(546,1644)(568,1646)(587,1647)
	(605,1648)(621,1649)(636,1649)
	(652,1650)(667,1649)(683,1649)
	(700,1648)(718,1647)(740,1646)
	(763,1644)(791,1642)(821,1639)
	(856,1635)(895,1630)(938,1624)
	(985,1617)(1034,1609)(1086,1599)
	(1139,1587)(1198,1572)(1254,1555)
	(1305,1539)(1349,1523)(1387,1509)
	(1419,1496)(1444,1484)(1464,1474)
	(1480,1465)(1492,1457)(1501,1450)
	(1508,1443)(1514,1437)(1520,1430)
	(1526,1423)(1534,1415)(1544,1406)
	(1558,1394)(1575,1381)(1596,1365)
	(1622,1347)(1653,1325)(1689,1300)
	(1728,1273)(1771,1243)(1814,1212)
	(1860,1178)(1902,1145)(1940,1115)
	(1972,1088)(1999,1065)(2022,1046)
	(2039,1030)(2053,1018)(2063,1008)
	(2071,1000)(2077,993)(2083,987)
	(2088,981)(2093,974)(2100,966)
	(2108,956)(2118,944)(2131,928)
	(2147,909)(2166,886)(2188,859)
	(2213,829)(2238,796)(2264,762)
	(2290,725)(2313,691)(2333,661)
	(2349,636)(2363,615)(2374,599)
	(2384,586)(2391,577)(2398,571)
	(2403,566)(2408,562)(2412,558)
	(2416,554)(2419,548)(2422,539)
	(2425,528)(2427,514)(2429,495)
	(2429,473)(2427,447)(2422,418)
	(2414,387)(2402,357)(2388,329)
	(2373,304)(2358,283)(2345,265)
	(2334,251)(2324,239)(2316,230)
	(2308,223)(2302,217)(2295,212)
	(2288,207)(2280,202)(2270,196)
	(2257,188)(2241,179)(2220,168)
	(2195,155)(2164,140)(2127,123)
	(2085,105)(2039,87)(1998,73)
	(1958,61)(1920,51)(1885,42)
	(1854,35)(1828,30)(1806,25)
	(1787,22)(1772,20)(1760,19)
	(1750,18)(1741,18)(1733,18)
	(1724,18)(1715,19)(1704,19)
	(1690,19)(1673,19)(1652,19)
	(1627,18)(1596,17)(1559,16)
	(1517,15)(1470,13)(1418,12)
	(1364,12)(1309,12)(1257,13)
	(1208,15)(1164,16)(1125,17)
	(1092,18)(1064,19)(1041,19)
	(1021,19)(1006,19)(992,19)
	(981,18)(970,18)(960,18)
	(949,18)(936,19)(921,20)
	(903,22)(881,25)(855,30)
	(825,35)(789,42)(750,51)
	(707,61)(661,73)(614,87)
	(565,103)(521,119)(481,134)
	(448,147)(420,157)(397,165)
	(380,171)(366,174)(356,176)
	(349,176)(343,175)(339,174)
	(335,174)(330,175)(325,177)
	(317,182)(308,190)(295,202)
	(280,219)(261,241)(240,269)
	(215,302)(190,342)(164,387)
	(143,428)(125,470)(108,511)
	(94,550)(81,587)(71,620)
	(62,650)(55,677)(49,701)
	(44,723)(40,742)(37,760)
	(35,777)(33,793)(31,810)
	(29,826)(28,845)(26,864)
	(24,886)(22,911)(20,938)
	(18,969)(16,1003)(14,1041)
	(13,1081)(12,1124)(12,1168)
	(14,1212)(17,1254)(22,1292)
	(27,1326)(32,1356)(36,1381)
	(39,1402)(42,1420)(43,1433)
	(44,1444)(44,1452)(43,1458)
	(42,1463)(40,1466)(39,1468)
	(38,1470)(38,1473)(39,1475)
	(41,1479)(45,1483)(52,1489)
	(62,1497)(74,1506)(91,1517)
	(112,1529)(137,1543)(166,1558)
	(201,1573)(239,1587)
\path(239,1587)(278,1599)(318,1609)
	(357,1617)(395,1624)(431,1630)
	(464,1635)(494,1639)(522,1642)
	(546,1644)(568,1646)(587,1647)
	(605,1648)(621,1649)(636,1649)
	(652,1650)(667,1649)(683,1649)
	(700,1648)(718,1647)(740,1646)
	(763,1644)(791,1642)(821,1639)
	(856,1635)(895,1630)(938,1624)
	(985,1617)(1034,1609)(1086,1599)
	(1139,1587)(1198,1572)(1254,1555)
	(1305,1539)(1349,1523)(1387,1509)
	(1419,1496)(1444,1484)(1464,1474)
	(1480,1465)(1492,1457)(1501,1450)
	(1508,1443)(1514,1437)(1520,1430)
	(1526,1423)(1534,1415)(1544,1406)
	(1558,1394)(1575,1381)(1596,1365)
	(1622,1347)(1653,1325)(1689,1300)
	(1728,1273)(1771,1243)(1814,1212)
	(1860,1178)(1902,1145)(1940,1115)
	(1972,1088)(1999,1065)(2022,1046)
	(2039,1030)(2053,1018)(2063,1008)
	(2071,1000)(2077,993)(2083,987)
	(2088,981)(2093,974)(2100,966)
	(2108,956)(2118,944)(2131,928)
	(2147,909)(2166,886)(2188,859)
	(2213,829)(2238,796)(2264,762)
	(2290,725)(2313,691)(2333,661)
	(2349,636)(2363,615)(2374,599)
	(2384,586)(2391,577)(2398,571)
	(2403,566)(2408,562)(2412,558)
	(2416,554)(2419,548)(2422,539)
	(2425,528)(2427,514)(2429,495)
	(2429,473)(2427,447)(2422,418)
	(2414,387)(2402,357)(2388,329)
	(2373,304)(2358,283)(2345,265)
	(2334,251)(2324,239)(2316,230)
	(2308,223)(2302,217)(2295,212)
	(2288,207)(2280,202)(2270,196)
	(2257,188)(2241,179)(2220,168)
	(2195,155)(2164,140)(2127,123)
	(2085,105)(2039,87)(1998,73)
	(1958,61)(1920,51)(1885,42)
	(1854,35)(1828,30)(1806,25)
	(1787,22)(1772,20)(1760,19)
	(1750,18)(1741,18)(1733,18)
	(1724,18)(1715,19)(1704,19)
	(1690,19)(1673,19)(1652,19)
	(1627,18)(1596,17)(1559,16)
	(1517,15)(1470,13)(1418,12)
	(1364,12)(1309,12)(1257,13)
	(1208,15)(1164,16)(1125,17)
	(1092,18)(1064,19)(1041,19)
	(1021,19)(1006,19)(992,19)
	(981,18)(970,18)(960,18)
	(949,18)(936,19)(921,20)
	(903,22)(881,25)(855,30)
	(825,35)(789,42)(750,51)
	(707,61)(661,73)(614,87)
	(565,103)(521,119)(481,134)
	(448,147)(420,157)(397,165)
	(380,171)(366,174)(356,176)
	(349,176)(343,175)(339,174)
	(335,174)(330,175)(325,177)
	(317,182)(308,190)(295,202)
	(280,219)(261,241)(240,269)
	(215,302)(190,342)(164,387)
	(143,428)(125,470)(108,511)
	(94,550)(81,587)(71,620)
	(62,650)(55,677)(49,701)
	(44,723)(40,742)(37,760)
	(35,777)(33,793)(31,810)
	(29,826)(28,845)(26,864)
	(24,886)(22,911)(20,938)
	(18,969)(16,1003)(14,1041)
	(13,1081)(12,1124)(12,1168)
	(14,1212)(17,1254)(22,1292)
	(27,1326)(32,1356)(36,1381)
	(39,1402)(42,1420)(43,1433)
	(44,1444)(44,1452)(43,1458)
	(42,1463)(40,1466)(39,1468)
	(38,1470)(38,1473)(39,1475)
	(41,1479)(45,1483)(52,1489)
	(62,1497)(74,1506)(91,1517)
	(112,1529)(137,1543)(166,1558)
	(201,1573)(239,1587)
\thicklines
\path(239,1587)(278,1599)(318,1609)
\path(614,1212)(1814,612)(614,612)(614,1212)
\end{picture}
}
\caption{Triangle and tangent cones}
\labell{fig:tangent-cones}
\end{figure}
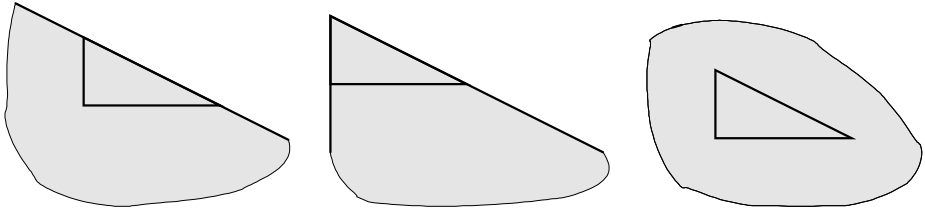

\begin{Example} \labell{ex:triangle}
Consider the triangle in $\R^2$ with vertices $(0,0)$, $(2,0)$,
and $(0,1)$.
The tangent cone at the hypotenuse is the (closed) half plane
consisting of all points lying below the line extending the hypotenuse;
the tangent cone at the top vertex consists of all rays subtended
from this vertex and pointing in the direction of the triangle;
if $F$ is the face consisting of the triangle itself,
then the tangent cone $\bfC_F$ is the whole plane.
See Figure \ref{fig:tangent-cones}.
\end{Example}

Let 
$ \sigma_1, \ldots, \sigma_d $
denote the facets (codimension one faces) of $\Delta$. (Warning:
Cappell and Shaneson use the symbols $\sigma_i$ to denote
the \emph{dual} objects to the facets, namely,
the one dimensional cones in the corresponding fan.)
Assume that $\Delta$ is simple, so that exactly $n$ facets intersect
at each vertex.  Let $\VertDelta$ denote the set of vertices of $\Delta$.
For each vertex $v \in \VertDelta$, let
$$ I_v \subset \{ 1, \ldots, d \} $$
encode the set of facets that meet at $v$, so that
$$ i \in I_v \quad \text{if and only if} \quad v \in \sigma_i . $$
Let $\alpha_{i,v}$, for $i\in I_v$, be edge vectors emanating from $v$;
concretely, assume that $\alpha_{i,v}$ lies along the unique edge at $v$ 
which is not contained in the facet $\sigma_i$.
(At the moment, the $\alpha_{i,v}$ are only determined
up to positive scalars.) 
In terms of the edge vectors, the tangent cone at a vertex $v$ is
$$ C_v = \{ v + \sum_{j \in I_v} x_j \alpha_{j,v} \ | \ x_j \geq 0
\text{ for all } j \} . $$

The polar decomposition theorem 
relates the characteristic function of the polytope to the
characteristic functions of convex polyhedral cones.  
As in the one dimensional case
\eqref{Lawrence for interval}, we cannot just consider the 
tangent cones, but we must make two modifications.
First, we must ``polarize" the tangent cones by flipping some of the
edge vectors.  Second, we must remove some facets.

To carry this out, we choose a vector $\xi \in V^*$
such that the pairings $\l< \xi , \alpha_{i,v} \r> $ are all non-zero;
we call it a ``polarizing vector"
and think of it as defining the ``upward" direction in $V$.
We ``polarize" the edge vectors so that they all point ``down":
for each vertex $v$ of $\Delta$
and each edge vector $\valpha_{i,v}$ emanating from $v$,
we define the corresponding \emph{polarized edge vector} to be
\begin{equation} \labell{polarization}
 \valpha^{\#}_{i,v} = \begin{cases}
\valpha_{i,v} & \text{ if } \l< \xi , \valpha_{i,v} \r> < 0, \\
- \valpha_{i,v} & \text{ if } \l< \xi , \valpha_{i,v} \r> > 0 .
\end{cases}
\end{equation}
Let
$$ \varphi_v = \{ j \in I_v \mid \l< \xi , \valpha_{j,v} \r> > 0 \}$$
denote the set of ``upward" edge vectors emanating from $v$, that is,
those edge vectors that get flipped in the polarization process
\eqref{polarization}.
The \emph{polarized tangent cone} to $\Delta$ at $v$
is obtained from the tangent cone $\bfC_v$ by flipping the $j$th edge
and removing the $j$th facet for each $j \in \varphi_v$:
\begin{equation} \labell{Cvsharp}
\bfC_v^\# =
\left\{ \left. v + \sum_{j \in I_v} x_j \alpha_{j,v}^\# \ \right| \ 
\begin{array}{l}
x_j \geq 0 \ \text{ if } j \in I_v \ssminus \varphi_v , \ \text{ and } \\
x_j   >  0 \ \text{ if } j \in \varphi_v \end{array}
\right\} .
\end{equation}

Recall that the characteristic function of a set $A \subset \R^n$ is
$$ \mathbf{1}_A(x) = \begin{cases} 1 & x \in A \\ 0 & x \not\in A .
\end{cases}$$

\begin{refTheorem}[Lawrence, Varchenko] \labell{thm:Lawrence}
\begin{equation}
\labell{Lawrence}
{\mathbf 1}_{\Delta}(x)
  = \sum_{v}(-1)^{|\varphi_v|}{\mathbf 1}_{\bfC_v^\#} (x).
\end{equation}
\end{refTheorem}

This decomposition was proved by Lawrence and Varchenko;
see \cite{lawrence,GV,V}.
A version for non-simple polytopes appeared in \cite{Haase}. 

In preparation for our proof of Theorem \ref{thm:Lawrence}
we introduce some notation.  

Let $E_1,\ldots,E_N$ be all the different hyperspaces in $V^*$
that are perpendicular to edges of $\Delta$ under the pairing
between $V$ and $V^*$.  That is,
\begin{equation} \labell{walls}
\{ E_i \ | \ 1 \leq i \leq N \} = \{ \ker \alpha_{j,v} \ | \ 
 v \in \VertDelta, \ j \in I_v \}.
\end{equation}
(For instance, if no two edges of $\Delta$ are parallel,
then the number $N$ of such hyperplanes
is equal to the number of edges of $\Delta$.)
A vector $\xi$ can be taken to be a ``polarizing vector" if and 
only if it belongs to the complement
\begin{equation} \labell{tDelta}
V^*_\Delta = V^* \ssminus \left( E_1 \cup \ldots \cup E_N \right).
\end{equation}
The connected components of this complement are called \emph{chambers}.
The signs of the pairings $\l< \xi , \alpha_{i,v} \r> $ 
only depend on the chamber containing $\xi$.

\begin{Remark} \labell{rk:chambers}
For a Hamiltonian action of a torus $T$ on a symplectic manifold $M$,
one similarly obtains chambers in the Lie algebra $\t$ of $T$
from the isotropy weights
$\alpha_{j,p}$ at the fixed points for the action.
When $M$ is a toric variety corresponding to the polytope $\Delta$,
this gives the same notion of chambers
as we have just described.  For a flag manifold $M \cong G/T$,
where $G$ is a compact Lie group and $T$ is a maximal torus,
the isotropy weights $\alpha_{j,p}$
are the roots of $G$, and the corresponding chambers are precisely the
interiors of the Weyl chambers.
\end{Remark}

Now suppose that $\xi$ belongs to exactly
one of the ``walls" in \eqref{walls}.
Let $e$ be an edge of $\Delta$ that is perpendicular to this wall.
Let
$$ I_e \subset \{ 1 , \ldots , d \} $$
correspond to the facets whose intersection is $e$, so that
$$ i \in I_e \quad \text{ if and only if } \quad e \subset \sigma_i .$$

Let $v$ be an endpoint of $e$.
The edge vectors at $v$ are $\alpha_{j,v}$,
for $j \in I_e$, and an edge vector that lies along $e$,
which we denote $\alpha_{e,v}$.  Note that
$$ \left< \xi , \alpha_{e,v} \right> = 0 . $$
Define
$$ \varphi_e = \{ j \in I_e \ | \ \left< \xi , \alpha_{j,v} \right> > 0 \}$$
and, for each $j \in I_e$,
$$ \alpha_{j,v}^\# = \begin{cases}
  \alpha_{j,v} & \text{ if } \left< \xi , \alpha_{j,v} \right> < 0 \\
 -\alpha_{j,v} & \text{ if } \left< \xi , \alpha_{j,v} \right> > 0.
\end{cases}$$
Define the \emph{polarized tangent cone} to $\Delta$ at the edge $e$ 
to be
$$ \bfC_e^\# = \left\{ \left. v + x_e \alpha_e 
   + \sum_{j \in I_e} x_j \alpha_{j,v}^\# \ \right| \ \
   \begin{array}{l}
   x_e \in \R,\\  x_j \geq 0 \ \text{ if } j \in I_e \ssminus \varphi_e, 
   \ \text{ and }\\ 
                  x_j > 0 \ \text{ if } j \in \varphi_e.\end{array} \right\} $$
Let $v'$ be the other endpoint of $e$.  We can normalize
the edge vectors such that $\alpha_{e,v'} = - \alpha_{e,v}$
and $\alpha_{j,v} - \alpha_{j,v'} \in \R \alpha_{e,v}$.
The set $\varphi_e$ and the cone $C_e^\#$ are independent 
of the choice of endpoint $v$ of $e$.

\begin{proof}[Proof of Theorem \ref{thm:Lawrence}]
Pick any polarizing vector $\xi \in V^*_\Delta$. Let $v \in \VertDelta$
be the vertex for which $\left< \xi , v \right> $ is maximal.
Then none of the $\alpha_{j,v}$'s are flipped,
and so $\bfC_v^\sharp = \bfC_v$.
For any other vertex $u \in \VertDelta$, at least one of the $\alpha_{j,u}$'s
is flipped, and so $\bfC_u^\sharp \cap \bfC_u = \emptyset $.
So the polytope $\Delta$
is contained in the polarized tangent cone $\bfC_v^\#$ at $v$
and is disjoint from the polarized tangent cone $\bfC_u^\#$ 
for all other $u \in \VertDelta$, and the equation \eqref{Lawrence}, 
when evaluated at $x \in \Delta$, reads $1=1$.

Suppose now that $x \not \in \Delta$.
The set of vectors $\xi$ which separate $x$ from $\Delta$,
that is, such that
\begin{equation} \labell{separate}
 \left< \xi , x \right> > \max\limits_{y \in \Delta}
 \left< \xi , y \right>, 
\end{equation}
is open in $V^*$.  Choose a polarizing vector $\xi \in V_\Delta^*$
that satisfies \eqref{separate}.
Then $x$ is not in the polarized tangent cone $\bfC_v^\#$
for any $v \in \VertDelta$.  Equation \eqref{Lawrence} for the 
polarizing vector $\xi$, when evaluated at $x$, reads $0=0$.

We finish by showing that, when the polarizing vector $\xi$ crosses 
a single wall $E_j$ in $V^*$, the right hand side of \eqref{Lawrence} 
does not change.

If $E_j$ is not perpendicular to any of the edge vectors at $v$, 
the signs of $\langle \xi , \alpha_{j,v} \rangle $ do not change,
so the polarized tangent cone $C_v^\#$ does not change
as $\xi$ crosses the wall.
The vertices whose contributions to the right hand side of
\eqref{Lawrence} change as $\xi$ crosses $E_j$ come in pairs,
because each edge of $\Delta$ that is perpendicular to $E_j$
has exactly two endpoints.

For each such vertex $v$, denote by $\bfS_v(x)$ and $\bfS_v'(x)$
its contributions to the right hand side of \eqref{Lawrence}
before and after $\xi$ crossed $E_j$.
Let $e$ be an edge perpendicular to $E_j$ and $v$ an endpoint of $e$.
Let $\bfS_e^\#(x)$ be the characteristic function of the
polarized tangent cone $C_e^\#$ corresponding to the value of $\xi$
as it crosses $E_j$.
The difference $\bfS_v(x)-\bfS_v'(x)$ is plus/minus $\bfS_e^\#(x)$.
For the \emph{other} endpoint $v$ of $e$,
the difference $\bfS_v(x)-\bfS_v'(x)$ is minus/plus $\bfS_e^\#(x)$.
So the differences $\bfS_v(x)-\bfS_v'(x)$,
for the two endpoints $v$ of $e$,  sum to zero.
\end{proof}

\begin{Remark} \labell{polar measure}
If we multiply both sides of (\ref{Lawrence}) by Lebesgue
measure  $dx$ we obtain a 
formula for ${\bf 1}_\Delta dx$ supported on $\Delta$ in terms of  
an alternating sum of the measures ${\mathbf 1}_{\bfC_v^\#} (x)dx$.
This formula (which is a special case of ``Filliman duality" 
\cite{filliman}) allows us to 
express the integral of any compactly supported continuous 
function $f$ over the polytope
in terms of its integrals over the cone $\bfC_v^\#$. 
From the point of view of measure theory, 
the missing facets of $\bfC_v^\#$ are irrelevant as they have
measure zero; what is important is the  change in the
direction of some of the edges at a vertex and the sign $(-1)^{|\varphi_v|}$ 
associated to the vertex. 
\end{Remark}

\begin{Remark}
If the $\Delta$ is a polytope with non-singular fan,
then (up to factors of $2\pi$) the measure ${\bf 1}_\Delta dx$ 
is the Duistermaat-Heckman measure of the associated toric manifold, 
and the vertices of 
$\Delta$ are the images of the fixed points of the torus 
acting on this manifold. In this case, the equality \eqref{Lawrence}
multiplied by Lebesgue measure $dx$ becomes a special case of 
the Guillemin-Lerman-Sternberg (G-L-S) formula.
The G-L-S formula expresses the Duistermaat-Heckman measure for a
Hamiltonian torus action in terms of an alternating sum 
of the Duistermaat-Heckman measure
associated to the linearized action along the components 
of the fixed point set. See \cite[section 3.3]{GLS}. This formula, in turn, 
can be deduced from the fact that a Hamiltonian torus action 
is cobordant (in an appropriate sense) to a union of its 
linearized actions along the components of its 
fixed point set.  See \cite[Chapter 4]{GGK}. 
The role of the polarizing vector $\xi\in V^*$ is 
played by a vector $\eta$ in the Lie algebra $\bf t$ such that 
the $\eta$-component $\Phi^\eta$ of the moment map $\Phi$
 is proper and bounded from below.  
\end{Remark}

\section{Sums and integrals over simple polytopes}
\labell{sec:first formula}

Formulas for sums and integrals of exponential functions
over simple polytopes appeared in \cite{BV,BV:partition}.
In this section we deduce such formulas from the ``polar decomposition".

We work in a vector space $V$ with a lattice $V_\Z$.
One may identify $V$ with $\R^n$ and $V_\Z$ with $\Z^n$, but we prefer
to use notation that is independent of the choice of a basis.
The \emph{dual} lattice is
\begin{equation} \labell{dual lattice} 
V^*_\Z = \{ u \in V^* \ | \ \left< u , \alpha \right> \in \Z
\, \text{ for all } \, \alpha \in V_\Z \}.
\end{equation}

\begin{Remark} \labell{chambers}
When our polytope is viewed as associated to a toric variety,
the vector space $V$ is the dual $\t^*$ of the Lie algebra $\t$ 
of a torus $T$, the lattice $V_\Z$ is the weight lattice in $\t^*$,
and $V_\Z^*$ is the kernel of the exponential map $\t \to T$.
\end{Remark}


Let $\bfC_v$ be a simple convex polyhedral cone in $V$, that is,
a set of the form
\begin{equation} \labell{Cv}
 \bfC_v = \{ v + \sum_{j=1}^n x_j \alpha_j \ | \ x_j \geq 0 
   \text{ for all } j \} 
\end{equation}
where $\alpha_1,\ldots,\alpha_n$ are a basis of $V$. 
Equivalently, we can write
\begin{equation} \labell{Cv as intersection}
 \bfC_v = 
 \{ x \ | \ \left< u_i , x \right> + \lambda_j \geq 0 \quad , \quad 
 j=1,\ldots,n \},
\end{equation}
where $u_1,\ldots,u_n$ are a basis of $V^*$.
We can pass from one description to another by setting $u_1,\ldots,u_n$
to be the dual basis to $\alpha_1,\ldots,\alpha_n$ and vice versa.
Geometrically, the vectors $u_i \in V^*$ are the inward normal vectors
to the facets of $\bfC_v$, and they encode the \emph{slopes} of the facets;
the real numbers $\lambda_i$ then determine the \emph{locations}
of the facets; the ``edge vectors" $\alpha_i$ generate the edges
of $\bfC_v$. 

A priori, the $u_i$'s and the $\alpha_i$'s are only determined 
up to multiplication by positive scalars.
To fix a particular normalization, we assume that the cone $\bfC_v$ 
is \emph{rational}, that is, that the $\alpha_i$'s can be chosen to be 
elements of the lattice $V_\Z$, or, equivalently, the $u_i$'s 
can be chosen to be elements of $V_\Z^*$.
We choose the normal vectors $u_i$ to be primitive lattice elements,
which means that each $u_i$ is in $V_\Z^*$ and is not equal to the product 
of an element of $V_\Z^*$ by an integer greater than one.
We choose $\alpha_1,\ldots,\alpha_n$ to be the dual basis 
to $u_1,\ldots,u_n$.

Although the $u_i$'s are primitive lattice elements, they may generate
a lattice that is coarser than $V_\Z^*$. 
The $\alpha_i$'s then generate a lattice that is \emph{finer} than $V_\Z$;
in particular, the $\alpha_i$'s themselves might not be in $V_\Z$.
See Figure \ref{fig:nonreg-orthant}.
The cone $\bfC_v$ is \emph{non-singular} if the $u_i$'s generate the lattice 
$V_\Z^*$, or, equivalently, the $\alpha_i$'s generate $V_\Z$. 
\begin{figure}[ht]
\setlength{\unitlength}{0.00083333in}
\begingroup\makeatletter\ifx\SetFigFont\undefined%
\gdef\SetFigFont#1#2#3#4#5{%
  \reset@font\fontsize{#1}{#2pt}%
  \fontfamily{#3}\fontseries{#4}\fontshape{#5}%
  \selectfont}%
\fi\endgroup%
{\renewcommand{\dashlinestretch}{30}
\begin{picture}(3380,2337)(0,-10)

\texture{0 0 0 888888 88000000 0 0 80808
        8000000 0 0 888888 88000000 0 0 80808
        8000000 0 0 888888 88000000 0 0 80808
        8000000 0 0 888888 88000000 0 0 80808 }
\thinlines
\shade\path(3368,2112)(1268,2112)(218,12)(2568,112)(3168,912)(3368,2112)

\put(68,2112){\blacken\ellipse{80}{80}}
\put(68,1512){\blacken\ellipse{80}{80}}
\put(68,912){\blacken\ellipse{80}{80}}
\put(68,312){\blacken\ellipse{80}{80}}

\put(668,2112){\blacken\ellipse{80}{80}}
\put(668,1512){\blacken\ellipse{80}{80}}
\put(668,912){\blacken\ellipse{80}{80}}
\put(668,312){\blacken\ellipse{80}{80}}

\put(1268,2112){\blacken\ellipse{80}{80}}
\put(1268,1512){\blacken\ellipse{80}{80}}
\put(1268,912){\blacken\ellipse{80}{80}}
\put(1268,312){\blacken\ellipse{80}{80}}

\put(1868,2112){\blacken\ellipse{80}{80}}
\put(1868,1512){\blacken\ellipse{80}{80}}
\put(1868,912){\blacken\ellipse{80}{80}}
\put(1868,312){\blacken\ellipse{80}{80}}

\put(2468,2112){\blacken\ellipse{80}{80}}
\put(2468,1512){\blacken\ellipse{80}{80}}
\put(2468,912){\blacken\ellipse{80}{80}}
\put(2468,297){\blacken\ellipse{80}{80}}

\put(3068,2112){\blacken\ellipse{80}{80}}
\put(3068,1512){\blacken\ellipse{80}{80}}
\put(3068,912){\blacken\ellipse{80}{80}}
\put(3068,312){\blacken\ellipse{80}{80}}

\put(2218,1787){\makebox(0,0)[lb]{\smash{{{\SetFigFont{12}{14.4}{\rmdefault}{\mddefault}{\updefault}$\sss u_1$}}}}}
\put(2218,1787){\makebox(0,0)[lb]{\smash{{{\SetFigFont{12}{14.4}{\rmdefault}{\mddefault}{\updefault}$\sss u_1$}}}}}
\put(1268,637){\makebox(0,0)[lb]{\smash{{{\SetFigFont{12}{14.4}{\rmdefault}{\mddefault}{\updefault}$\sss u_2$}}}}}
\put(923,1827){\makebox(0,0)[lb]{\smash{{{\SetFigFont{12}{14.4}{\rmdefault}{\mddefault}{\updefault}$\sss \alpha_1$}}}}}
\put(1343,2207){\makebox(0,0)[lb]{\smash{{{\SetFigFont{12}{14.4}{\rmdefault}{\mddefault}{\updefault}$\sss \alpha_2$}}}}}

\thicklines
\path(3368,2112)(1268,2112)(218,12)

\path(2468,2112)(2468,1512)
\path(2438.000,1632.000)(2468.000,1512.000)(2498.000,1632.000)

\path(1268,2112)(968,1512)
\path(994.833,1632.748)(968.000,1512.000)(1048.498,1605.915)

\path(1268,2112)(1568,2112)
\path(1448.000,2082.000)(1568.000,2112.000)(1448.000,2142.000)

\path(668,912)(1868,312)
\path(1747.252,338.833)(1868.000,312.000)(1774.085,392.498)

%
%
%

\end{picture}
}
\caption{}
\labell{fig:nonreg-orthant}
\end{figure}

Let $dx$ denote Lebesgue measure on $V$, normalized so that the measure
of $V/V_\Z$ is one.  Consider an exponential function $f \colon V \to \C$
given by
\begin{equation} \labell{what is f}
 f(x) = e^{\left< \xi , x \right> } ,
\end{equation}
where $\xi \in V^*_\C$ is such that
$$ \Re \left< \xi , \alpha_i \right>  < 0 \quad \text{for all $i$}.$$
Then the integral of $f$ over the cone $\bfC_v$ and the sum of $f$ 
over the lattice points in $\bfC_v$ both converge.

If the cone $\bfC_v$ is non-singular and $v$ is in the lattice $V_\Z$
then the map
$$ (t_1,\ldots,t_n) \mapsto v + \sum t_j \alpha_j $$
sends the standard positive orthant 
$$ \R_+^n = \{ (t_1,\ldots,t_n) \ | \ t_j \geq 0 \text{ for all } j \} $$
onto $\bfC_v$ and sends $\R_+^n \cap \Z^n$ onto $\bfC_v \cap V_\Z$.
In particular, it takes the standard Lebesgue measure on $\R^n$ to 
the measure $dx$ on $V$. So
\begin{equation} \labell{dot}
\begin{aligned}
\int_{\bfC_v} f(x) dx &= \int_0^\infty \cdots \int_0^\infty 
  e^{ \left< \xi , v + \sum t_j \alpha_j \right> } dt_1 \cdots dt_n \\
 &= e^{\left< \xi , v \right> } 
   \prod_{j=1}^n \int_0^\infty e^{ t \left< \xi , \alpha_j \right> } dt_j \\
 &= e^{\left< \xi,v \right> } 
   \prod_{j=1}^n - \frac{1}{\left< \xi , \alpha_j\right> } 
\end{aligned}
\end{equation}
and
\begin{equation} \labell{dot2} 
\begin{aligned}
\sum_{\bfC_v \cap V_\Z} f &= \sum_{k_1=0}^\infty \cdots \sum_{k_n=0}^\infty 
   e^{\left< \xi , v+\sum k_j \alpha_j \right> } \\
 &= e^{ \left< \xi , v \right> } \prod_{j=1}^n \sum_{k_j=0}^\infty 
    \left( e^{ \left< \xi , \alpha_j \right> } \right)^{k_j} \\
 &= e^{ \left< \xi , v \right> } \prod_{j=1}^n 
   \frac{1}{1 - e^{\left< \xi , \alpha_j \right> } }.
\end{aligned}
\end{equation}

A crucial ingredient in extending the Khovanskii-Pukhlikov formula 
to the case of simple polytopes is an extension of the formulas
\eqref{dot} and \eqref{dot2}
to the case that the cone $\bfC_v$ is not non-singular.

We associate to $\bfC_v$ the finite abelian group 
$$ \Gamma = V_\Z^* / \span_\Z \{ u_i \} .$$
Note that the group $\Gamma$ is trivial if and only if the cone $\bfC_v$ 
is non-singular.
Also note that $e^{2\pi i \left< \gamma, x \right> } $
is well defined whenever $\gamma \in \Gamma$
and $x \in \span_\Z \{ \alpha_i \} $.

\begin{Remark*}
The toric variety associated to the cone is $\C^n/\Gamma$, 
where $\Gamma$ acts on $\C^n$ through the homomorphism $\Gamma \to (S^1)^d$
given by $\gamma \mapsto \left( 
        e^{2\pi i \left< \gamma , \alpha_1 \right> },
\ldots, e^{2\pi i \left< \gamma , \alpha_n \right> } \right) $.
\end{Remark*}

We have the following generalization of Formula \eqref{dot2}.
Suppose that the vertex of $C_v$ satisfies $v \in \span_\Z \{ \alpha_i \}$.
Then
\begin{equation} \labell{claim1}
\sum_{x \in \bfC_v \cap V_\Z} e^{\l< \xi , x \r> }
\; = \; e^{\l< \xi,v \r> } \cdot \frac{1}{|\Gamma|} \sum_{\gamma \in \Gamma}
e^{2\pi i \l< \gamma,v \r> } \prod_{j=1}^n
\frac{1} 
     {1 - e^{2\pi i \l< \gamma,\valpha_j \r> } e^{\l< \xi,\valpha_j \r> } }.
\end{equation}

\begin{proof}
The main step is to transform the left hand side of \eqref{claim1}
into a summation over elements of the finer lattice
$\span_\Z \{ \alpha_j \}$.
For each $x \in \span_\Z \{ \alpha_i \}$,
\begin{equation} \labell{characterX}
\gamma \mapsto e^{2\pi i \l< \gamma , x \r> }
\end{equation}
is a homomorphism from $\Gamma$ to $S^1$, 
and it is trivial if and only if $x \in V_\Z$.
Frobenius's theorem asserts that, for a finite group $\Gamma$,
the sum of the values of a non-trivial homomorphism $\Gamma \to S^1$
is zero.  It follows that, for $x \in \span_\Z \{ \alpha_i \}$,
\begin{equation} \labell{frobenius}
\frac{1}{|\Gamma|} \sum_{\gamma \in \Gamma} e^{2\pi i \l< \gamma,x \r> }
\; = \; \begin{cases}
 1 & \text{ if } x \in V_\Z ; \\
 0 & \text{ otherwise.}
\end{cases} 
\end{equation}
By \eqref{frobenius}, the left hand side of \eqref{claim1} is equal to
\begin{equation} \labell{almost there}
\sum_{x \in \bfC_v \cap \span_\Z \{ \alpha_j \} } e^{\l< \xi, x \r> }
\cdot \frac{1}{|\Gamma|} \sum_{\gamma \in \Gamma} 
e^{2\pi i \l< \gamma , x \r> } \ .
\end{equation}
Writing
$$ x = v + \sum k_j \valpha_j , $$
this becomes
$$
e^{ \l< \xi, v \r> } \cdot \frac{1}{|\Gamma|} \sum_{\gamma \in \Gamma}
e^{2\pi i \l< \gamma , v \r> }
\prod_{j=1}^n \sum_{k=0}^\infty 
   e^{2\pi i k \l< \gamma , \valpha_j \r> } e^{k \l< \xi , \valpha_j \r> } ,
$$
which is equal to the right hand side of \eqref{claim1} 
by the formula for the sum of a geometric series.
\end{proof}

We also have the following generalization of Formula \eqref{dot}.
\begin{equation} \labell{claim2}
\int_{\bfC_v} e^{\l< \xi , x \r> }  dx
\; = \; e^{\l< \xi,v \r> } \cdot \frac{1}{|\Gamma|} \prod_{j=1}^n
- \frac{1}{\l< \xi,\valpha_j \r> } .
\end{equation}

\begin{proof}
We perform the change of variable $x = v + \sum_{j=1}^n t_j \valpha_j$.
Then $x \in \bfC_v$ if and only if $t = (t_1,\ldots,t_n)$ 
belongs to the positive orthant $\R_+^n$.
The inverse transformation is
$$ t_j = \left< u_j , x-v \right> ; $$
its Jacobian is $ [ V_\Z^* : \span_\Z u_j ] = | \Gamma | $. So
\begin{eqnarray*}
\int_{\bfC_v} e^{\l< \xi , x \r> } dx & = &
\frac{1}{|\Gamma|} \int_{\R_+^n} e^{\left< \xi , v+\sum t_j \alpha_j \right>} 
   dt_1 \cdots dt_n \\
 & = & e^{\l<  \xi , v \r> } \cdot \frac{1}{|\Gamma|}
\prod_{j=1}^n \int_0^\infty e^{t \l< \xi , \valpha_j \r> } dt \\
 & = & e^{\l<  \xi , v \r> } \cdot \frac{1}{|\Gamma|}
\prod_{j=1}^n - \frac{1}{\l< \xi,\valpha_j \r> } .
\end{eqnarray*}
\end{proof}

To apply Formulas \eqref{claim1} and \eqref{claim2} to the ``polar
decomposition" of a polytope, we need to consider ``polarized cones".
Suppose that $\xi \in V^*_\C$ satisfies
$\Re \left< \xi , \alpha_j \right> \neq 0$ for all $j$.  
As in section \ref{sec:decompose},
let
$$ \varphi_v = \{ j \ | \ \Re \left< \xi , \alpha_j \right> > 0 \} ;$$
for each $j$, let
\begin{equation} \labell{alpha j sharp}
\alpha_j^\sharp = \begin{cases} 
   \alpha_j & \text{ if } j \not\in \varphi_v \\
 - \alpha_j & \text{ if } j     \in \varphi_v ;
\end{cases}
\end{equation}
and let
$$ \bfC_v^\sharp = \left\{ \left. 
                   v + \sum_{j=1}^n x_j \alpha_j^\sharp \ \right| \ 
\begin{array}{l}
x_j \geq 0 \text{ if } j \not\in \varphi_v, \text{ and } \\
x_j > 0 \text{ if } j \in\varphi_v.
\end{array} \right\}.$$
Then the integral of $f$ over $\bfC_v^\sharp$ and the sum of $f$
over the lattice points in $\bfC_v^\sharp$ converge.  We compute
this integral and this sum:

\begin{equation} \labell{BB1}
\begin{aligned}
\int_{\bfC_v^\sharp} e^{\left< \xi , x \right> } dx
 & = e^{\left< \xi , v \right> } \cdot \frac{1}{|\Gamma|} 
     \prod_{j=1}^n - \frac{1}{\langle \xi , \alpha^\#_j \rangle }
 & \quad \text{ by \eqref{claim2} } \\
 & = (-1)^{|\varphi_v|} e^{\left< \xi , v\right> } \cdot \frac{1}{|\Gamma|}
     \prod_{j=1}^n - \frac{1}{\left< \xi , \alpha_j \right> } 
 & \quad \text{by \eqref{alpha j sharp}.}
\end{aligned}
\end{equation}
Because $v \in \span_\Z \{ \alpha_j \} $,
$$\begin{aligned}\
\bfC_v^\sharp \cap ( \span_\Z \{ \alpha_j \} )
 &= \left\{ \left. v + \sum m_j \alpha_j^\sharp \ \right| \ 
{\begin{array}{l} m_j \in \Z, \\
                  m_j \geq 0 \text{ if } j \not\in \varphi_v, \\
                  m_j > 0    \text{ if } j     \in \varphi_v 
\end{array}} \right\} \\
 &= \left\{ \left. v_\shift + \sum m_j \alpha_j^\sharp \ \right| \ 
 {\begin{array}{l}
 m_j \in \Z, \\ 
 m_j \geq 0 \text{ for all } j 
 \end{array}} \right\}
\end{aligned} $$
where
\begin{equation} \labell{star5}
 v_\shift = v + \sum_{j \in \varphi_v} \alpha_j^\sharp .
\end{equation}
So
\begin{equation} \labell{same integral}
 \bfC_v^\sharp \cap V_\Z = \ol{\bfC}^\sharp_{v,\shift} \cap V_\Z 
\end{equation}
where
$$ \ol{\bfC}^\sharp_{v,\shift} = \{ v_\shift + \sum x_j \alpha_j^\sharp 
    \ | \ x_j \geq 0 \text{ for all } j   \} .$$
We have
\begin{equation} \labell{BB2}
\begin{aligned}
\sum_{x \in {\bfC}_v^\sharp \cap V_\Z} e^{\left< \xi , x \right> }
 &= \sum_{x \in \ol{\bfC}^\sharp_{v,\shift} \cap V_\Z} 
                                        e^{\left< \xi , x \right>} 
 \qquad \text{ by \eqref{same integral} } \\
 &= e^{\left< \xi , v_\shift \right>} \cdot \frac{1}{| \Gamma |} 
    \sum_{\gamma \in \Gamma} e^{ 2 \pi i {\left< \gamma , v_\shift \right> } }
    \prod_{j=1}^n 
    \frac{1}{1 - e^{ 2 \pi i \left< \gamma , \alpha_j^\sharp \right> }
                 e^{\left< \xi , \alpha_j^\sharp \right> } }
 \qquad \text{ by \eqref{claim1} } \\
 &= e^{\left< \xi , v \right> } \cdot\frac{1}{|\Gamma|} \sum_{\gamma \in \Gamma}
    e^{ 2 \pi i \left< \gamma , v \right> }
    \prod_{j \not\in \varphi_v } 
    \frac{1} {1 - e^{2 \pi i\left< \gamma , \alpha_j \right> } 
                  e^{\left< \xi , \alpha_j \right> } }
    \prod_{j \in \varphi_v } 
    \frac{e^{-2 \pi i \left< \gamma , \alpha_j \right> } 
          e^{-\left< \xi , \alpha_j \right> } }
         {1 - e^{-2\pi i \left< \gamma,\alpha_j \right> } 
              e^{-\left< \xi , \alpha_j \right> } }
 \qquad \text{by \eqref{star5} and \eqref{alpha j sharp}}  \\
 &= (-1)^{|\varphi_v|} e^{\left< \xi , v \right>} \cdot \frac{1}{|\Gamma|}
    \sum_{\gamma \in \Gamma} e^{2 \pi i \left< \gamma , v \right> }
    \prod \frac{1}{1 - e^{2 \pi i \left< \gamma , \alpha_j \right> } 
                       e^{ \left< \xi , \alpha_j \right> } } \\
 & \qquad \qquad  \text{ by applying the relation } 
   \frac{e^x}{1-e^x} = - \frac{1}{1-e^{-x}} \text{ to } x = 
    - 2 \pi i \left< \gamma,\alpha_j\right> - \left< \xi , \alpha_j \right> 
   \text{ for } j \in \varphi_v.
\end{aligned}
\end{equation}

We can now reproduce Brion-Vergne's formulas for simple polytopes.
Let $\Delta \subset V$ be a simple polytope.
Suppose that $\xi \in V^*_\C$ satisfies 
$\Re \left< \xi , \alpha_{j,v} \right> \neq 0 $
for all $v \in \Vert(\Delta)$ and all $j \in I_v$.
With the notation of Section \ref{sec:decompose},
\begin{equation} \labell{B1}
\begin{aligned}
\int_\Delta e^{\left< \xi , x \right> } dx 
  &= \sum_{v \in \Vert(\Delta)} (-1)^{|\varphi_v|} 
     \int_{\bfC_v^\#} e^{ \left< \xi , x \right> } dx
  & \quad \text{ by \eqref{Lawrence} } \\
  &= \sum_{v \in \Vert(\Delta)}
  e^{\left< \xi , v \right> } \cdot \frac{1}{|\Gamma_v|} \prod_{j \in I_v}
  - \frac{1}{\left< \xi , \alpha_{j,v} \right> }
   & \quad \text{ by \eqref{BB1} }.
\end{aligned}
\end{equation}
Similar formulas appeared in \cite[Proposition 3.10]{BV}
and \cite[p.801, Theorem, part (ii)]{BV:partition}.
\comment{similar or same?}

Also,
\begin{equation} \labell{B2}
\begin{aligned}
  \sum_{x \in \Delta \cap V_\Z } e^{\left< \xi , x \right> } 
   &= \sum_{v \in \Vert(\Delta)} (-1)^{|\varphi_v|} 
      \sum_{x \in \bfC_v^\# \cap V_\Z} e^{ \left< \xi , x \right> }
   & \quad \text{ by \eqref{Lawrence}} \\ 
   &= \sum_{v \in \Vert(\Delta)} 
   e^{\left< \xi , v \right> } \cdot \frac{1}{|\Gamma_v|}
  \sum_{\gamma \in \Gamma_v} e^{2 \pi i \left< \gamma , v \right> }
  \prod_{j \in I_v} \frac{1}
        {1 - e^{2 \pi i \left< \gamma,\alpha_{j,v} \right> } 
             e^{\left< \xi , \alpha_{j,v} \right> } } 
  & \quad \text{ by \eqref{BB2}} \\
 &= \sum_{v \in \Vert(\Delta)} e^{\left< \xi , v \right> }
    \Td_v \left( \{ - \left< \xi , \alpha_{j,v} \right> \} \right)
    \cdot \frac{1}{|\Gamma_v|} \prod_{j \in I_v} 
                     - \frac{1}{\left< \xi , \alpha_{j,v} \right> }
\end{aligned}
\end{equation}
where
\begin{equation} \labell{Tdv is}
\Td_v\left( S \right) = 
 \sum_{\gamma \in \Gamma_v} e^{2 \pi i \left< \gamma,v \right>}
 \prod_{ j \in I_v } 
  \frac{S_j}{1 - e^{2 \pi i \left< \gamma,\alpha_{j,v}\right> } e^{-S_j} }
\quad \text{ for } \quad S = \{ S_j \}_{j \in I_v} .
\end{equation}
Similar formulas appeared in \cite[Proposition 3.9]{BV}
and \cite[p.801, Theorem, part (iii)]{BV:partition}


\begin{Remark} \labell{rk:cat}
We proved \eqref{B1} for $\xi$ outside the real hyperplanes
\begin{equation} \labell{cat}
\Re \l< \xi , \alpha_{j,v} \r> = 0 \quad,  \quad
   v \in \Vert(\Delta) \quad , \quad j \in I_v 
\end{equation}
in $V^*_\C$.  However, the left hand side of \eqref{B1} is analytic
for all $\xi \in V^*_\C$, and the right hand side is analytic outside 
the \emph{complex} hyperplanes
$$ \left< \xi , \alpha_{j,v} \right> = 0 .$$
By analytic continuation \eqref{B1} continues to hold for all $\xi$
outside these complex hyperplanes.  
Similarly, we proved \eqref{B2} for $\xi$ outside the real hyperplanes
\eqref{cat}, but by analytic continuation is remains true for all $\xi$
outside the \emph{complex} hyperplanes
\begin{equation} \labell{funny}
 \l< \xi , \alpha_{j,v} \r> = 2 \pi i \l< y,\alpha_{j,v} \r>
\quad , \quad y \in V^*_\Z .
\end{equation}
(Notice that these complex hyperplanes are contained in the real 
hyperplanes \eqref{cat}.)
We are grateful to A.\ Khovanskii 
for calling our attention to this approach. 
\end{Remark}

\begin{Remark} \labell{disk}
Note that the function $\Td_v(S)$ is analytic on the polydisk
$$ | S_j | < b_j \ , \ j \in I_v $$
in $\C^{I_v}$, where 
$$ b_j = \min_{\substack{y \in V_\Z^* \\ 
\left< y , \alpha_{j,v} \right> \neq 0}}
        \left\{ \, | 2 \pi \left< y , \alpha_{j,v} \right> | \, \right\}.$$
\end{Remark}

\section{Euler Maclaurin formulas for a simple polytope}
\labell{sec:EM simple polytope}

In this section we present Euler Maclaurin formulas for simple
lattice polytopes in arbitrary dimensions.  
As in the previous section, we work with an $n$ dimensional 
vector space $V$
with a lattice $V_\Z$.  Let $\Delta$ be a convex polytope in $V$ 
with $d$ facets, given by
\begin{equation} \labell{Delta}
\Delta = \bigcap_{i=1}^d \{ x \ | \ \l< x , u_i \r> + \lambda_i \geq 0 \}.
\end{equation}
The vectors $u_i \in V^*$ are inward normal vectors to the facets, 
and they encode the slopes of the facets; the real numbers $\lambda_i$ 
determine the location of the facets.
As before, we assume that the slopes of the facets are rational,
and we choose the normal vectors $u_i$ to be primitive elements
of the dual lattice $V_\Z^*$.  
We assume that the polytope $\Delta$ is simple, 
meaning that exactly $n$ facets intersect at each vertex.
We also assume that the $\lambda_i$'s are integers.

\begin{Remark}
If the vertices of $\Delta$ are lattice points,
the $\lambda_i$'s are integers.
If $\Delta$ has a non-singular fan and the $\lambda_i$'s are integers,
then the vertices of $\Delta$ are lattice points.
However, if $\Delta$ is simple but does not have a non-singular fan,
its vertices might not be lattice points even if the $\lambda_i$
are all integers.
\end{Remark}

For $h$ near $0$, the expanded polytope" 
\begin{equation} \labell{polytope}
\Delta(h) = \bigcap_{i=1}^d \{ x \in V \ \mid \
\l< x , u_i \r> + \lambda_i + h_i \geq 0 \}
\end{equation}
is obtained from $\Delta$ by shifting the half-spaces
defining its facets without changing their slopes.  

As before,
we normalize Lebesgue measure on $V$ so that a fundamental domain
with respect to the lattice $V_\Z$ has measure one. The integral
of a function $f$ over the expanded polytope $\Delta(h)$ 
is a function of $h_1,\ldots,h_d$.

For a polynomial $P$ in $d$ variables, the expression
\begin{equation} \labell{lambda expression}
\left. P \left( \deldel{h_1},\ldots,\deldel{h_d} \right) \right|_{h=0}
\int_{\Delta(h)} f
\end{equation}
makes sense if the integral is a smooth function of $h$
for $h$ near $0$.
If $P$ is a power series in $d$ variables, the expression
\eqref{lambda expression} makes sense if the resulting series converges.
If $P$ is an analytic function in $d$ variables, we interpret the
expression \eqref{lambda expression} by expanding $P$
into its Taylor series about the origin.

In what follows, the function $f$ can be taken to be
the product of a polynomial function with an exponential function
of the form $e^{ \left< \xi , x \right> }$ where $\xi \in V^*_\C$
is sufficiently small.

\smallskip\noindent\textbf{Khovanskii-Pukhlikov's formula 
                           for a polytope with a non-singular fan.}
%

The formula of Khovanskii and Pukhlikov
(see section 4 of \cite{KP2}),
translated to our notation, is the following formula.
\begin{equation} \labell{eq:KK}
 \sum_{\Delta \cap V_\Z} f  = \left.
 \Td \left( \deldel{h_1} \right)  \cdots \Td \left( \deldel{h_d} \right)
 \right|_{h=\lambda} \int_{\Delta(h)} f ,
\end{equation}
where the polytope $\Delta$ is integral and has a non-singular fan.

\smallskip\noindent
\textbf{Finite groups associated to the faces of a simple rational polytope.}
%

The facets of the polytope $\Delta$ are
$$ \sigma_i = \{ x \in \Delta \ | \ 
 \left< u_i , x \right> + \lambda_i = 0 \} \quad , \quad i=1,\ldots,d. $$
Because the polytope $\Delta$ is simple, each face $F$ of $\Delta$
can be uniquely described as an intersection of facets.  We let
$I_F \subset \{ 1 , \ldots , d \}$ denote the subset such that
$$ F = \bigcap\limits_{i \in I_F} \sigma_i .$$
The number of elements in $I_F$ is equal to the codimension of $F$.  
The relation
$$ F \mapsto I_F $$
is order (inclusion) reversing.  

For each vertex $v$ of $\Delta$, the vectors
$$ u_i \quad , \quad i \in I_v $$
form a basis of $V^*$. Let
$$ \alpha_{i,v} \quad , \quad i \in I_v $$
be the dual basis.
The $\alpha_{i,v}$'s are edge vectors at $v$, that is, they point
in the directions of the edges emanating from $v$.

The vector space $V_F$ normal to a face $F$ is the quotient of $V$
by $T_F = \{ r(x-y) \ | \ x,y \in F , r \in \R \}$.
Its dual is the subspace
\begin{equation} \labell{VFstar}
 V_F^* := \span \{ u_j \ | \ j \in I_F \}
\end{equation}
of $V^*$.  Let $\alpha_{j,F}$, $j \in I_F$, be the basis of $V_F$
that is dual to the basis $u_j$, $j \in I_F$, of $V_F^*$.

To each face $F$ of $\Delta$ we associate a finite abelian
group $\Gamma_F$ in the following way.  The lattice
$$ \span_\Z \{ u_i \ | \ i \in I_F \} \, \subset \, V_F^* \cap V_\Z^* $$
is a sublattice of $V_F^* \cap V_\Z^*$ of finite index.
The finite abelian group associated to the face $F$ is
the quotient
\begin{equation} \labell{GammaF}
\Gamma_F := 
\left( V_F^* \cap V_\Z^* \right) / \span_\Z \{ u_i \ | \ i \in I_F \} .
\end{equation}
For each $\gamma \in \Gamma_F$ and $j \in I_F$, 
the pairing $\l< \gamma , \valpha_{j,F} \r> $ 
is well defined modulo $1$, so 
$$ e^{2 \pi i\l< \gamma, \valpha_{j,F} \r> }$$ 
is well defined.

\begin{Remark}
The group $\Gamma_F$ measures the singularity of the toric variety 
associated to $\Delta$ along the stratum corresponding to $F$.
\end{Remark}

If $F \subseteq E$ are faces of $\Delta$, so that $I_E \subseteq I_F$, 
then $ \{ u_i \ | \ i \in I_E \} $ is a subset of
$\{ u_i \ | \ i \in I_F \}$.
Because these sets are bases of $V_E^*$ and $V_F^*$, we have
$$ V_E^* \subseteq V_F^* ,$$
and
$$ V_E^* \cap \span_\Z \{ u_i \ | \ i \in I_F \} 
   = \span_\Z \{ u_i \ | \ i \in I_E \} .$$
Hence, the natural map from  $\Gamma_E$ to $\Gamma_F$ is one to one,
and provides us with a natural inclusion map:
$$ \text{if \  $F \subseteq E$ \  then \  $\Gamma_E \subseteq \Gamma_F$.} $$

We define
\begin{equation} \labell{def:GammaFsharp}
\Gamma_F^\flat
   = \Gamma_F \ssminus
\bigcup_{\text{faces } E \text{ such that }  E \supsetneq F} \Gamma_E.
\end{equation}

Finally, note that, for each face $F$ and element $\gamma \in \Gamma_F$,
the number $e^{- 2 \pi i \left< \gamma , x \right> }$ is 
the same for all $x \in F$; we denote this number
$$ e^{- 2 \pi i \left< \gamma , F \right> }. $$

\smallskip\noindent
\textbf{Guillemin and Brion-Vergne formulas for a simple polytope.}
%

On any linear subspace $A$ of $V$ with rational slope
we normalize Lebesgue measure so that a fundamental domain 
with respect to the lattice $A \cap V_\Z$ has measure one.  
We shift this measure to any affine translate of $A$.
Integration over each face $F$ of $\Delta$ is defined
with respect to these measures.

For each face $F$ of $\Delta$, let
\begin{equation} \labell{F(h)}
 F(h) = \Delta(h) \cap \bigcap\limits_{i \in I_F}
 \{ x \ | \ \left< u_i , x \right> + \lambda_i + h_i = 0 \} 
\end{equation}
denote the corresponding face of the expanded polytope $\Delta(h)$.

Guillemin gives an Euler Maclaurin formula for a polytope
when the polytope
is expressed as the set of solutions of the equation
$k_1 \alpha_1 + \ldots + k_d \alpha_d = \mu$,
$k_1,\ldots,k_d \in \R_{\geq 0}$,
for some fixed integral vectors $\alpha_1,\ldots,\alpha_d,\mu$.
(See Theorem 1.3 and formula (3.28) of \cite{Gu}.)
When translated to our setup, his formula becomes the following formula.
\begin{equation} \labell{Guillemin}
\sum_{\Delta \cap V_\Z} f =
\sum_F \frac{1}{|\Gamma_F|} \sum_{\gamma \in \Gamma_F^\flat} 
\left.
e^{2 \pi i \left< \gamma , F \right> } 
\prod_{j \not \in I_F} \frac{\deldel{h_j}}{1 - e^{-\del/\del{h_j}}}
\prod_{j \in I_F}
\frac{1}{1 - e^{2\pi i \l< \gamma, \valpha_{j,F} \r> } e^{-\del/\del{h_j}}}
\right|_{h=0}
\int_{F(h)} f , 
\end{equation}
where the polytope $\Delta$ is simple and is given by \eqref{Delta}
where all the $\lambda_i$'s are integers.
\smallskip

Finally, the formula of Brion and Vergne in our notation is
\begin{equation} \labell{BV}
\sum_{\Delta \cap V_\Z} f =  
\left.  \sum_F \sum_{\gamma \in \Gamma_F^\flat }
e^{2 \pi i \left< \gamma , F \right> } 
\prod_{j \not \in I_F}
\frac{\deldel{h_j}}{1 - e^{- \del/\del{h_j}}}
\prod_{j \in I_F}
\frac{\deldel{h_j}}{1 - e^{2\pi i \l< \gamma,\valpha_{j,F} \r> }
e^{-\del/\del{h_j}}}
\right|_{h=0}
\int_{\Delta(h)} f .
\end{equation}
See \cite[Theorem 2.15]{BV} 
(where $\Delta$ is simple and integral), and, 
more generally, \cite{BV:partition}.

\begin{Remark} 
If the polytope $\Delta$ is integral
(which is a stronger requirement than the assumption that 
the $\lambda_i$'s be integers) then $e^{2\pi i \l< \gamma,F\r> } = 1$ 
for each face $F$ and each $\gamma \in \Gamma_F$.  
\end{Remark}


We now give a self-contained statement and an elementary proof
of the Guillemin-Brion-Vergne formulas.

\begin{Theorem} \labell{EM:exponentials}
Let $V$ be a vector space with a lattice $V_\Z$.
Let $V_\Z^* \subset V^*$ be the dual lattice. 
Let $\Delta \subset V$ be a simple rational polytope with $d$ facets.
Let $u_1,\ldots,u_d \in V^*_\Z$ be the primitive inward normals to the 
facets of $\Delta$. Let $\lambda_1,\ldots,\lambda_d$ be the real numbers
so that 
$$ \Delta = \bigcap_{i=1}^d \{ x \ | \ \left< u_i , x \right> + \lambda_i
   \geq 0 \}.$$
Suppose that $\lambda_1,\ldots,\lambda_d$ are integers.
For $h=(h_1,\ldots,h_d)$, let
$$ \Delta(h) = \bigcap_{i=1}^d \{ x \ | \ \l< u_i , x \r> 
   + \lambda_i + h_i \geq 0 \}.$$
For each face $F$ of $\Delta$, let $I_F \subset \{ 1, \ldots, d \}$
be the subset such that $F$ consists of those $x \in \Delta$
for which $\left< u_i , x \right> + \lambda_i = 0$ for all $i \in I_F$.
Let $\alpha_{i,F}$, for $i \in I_F$, be the basis of $V/T_F = V/\R(F-F)$
that is dual to the basis $u_i$, $i \in I_F$, of $V_F^* = T_F^0$. 
(In particular, if $v \in \Vert(\Delta)$ 
and $i \in I_v$ then $\alpha_{i,v}$ are the edge vectors emanating from $v$.) 
Let
$$\Gamma_F = (V_F^* \cap V_\Z^*) / \span_\Z \{ u_i \ | \ i \in I_F \}$$
be the finite group associated to the face $F$, and let
$\Gamma_F^\flat = \Gamma_F \ssminus \bigcup \Gamma_E$,
where the union is over all faces $E$ such that $E \supsetneq F$.
Let
$$ \Td_\Delta(S_1,\ldots,S_d) = \sum_F \sum_{\gamma \in \Gamma_F^\flat}
   e^{2 \pi i \left< \gamma , F \right> }
   \prod_{j \not\in I_F} \frac{S_j}{1 - e^{-S_j}} 
   \prod_{j \in I_F} \frac{S_j}
           {1 - e^{ 2 \pi i \l< \gamma , \alpha_{j,F} \r> } e^{-S_j} }.$$ 
Let $f \colon V \to \C$ be a quasi-polynomial function,
that is, a linear combination of functions of the form
$$ f(x) = p(x) e^{\left< \xi , x \right> } $$
where $p \colon V \to \C$ are polynomial functions
and where the exponents $\xi \in V_\C^*$ satisfy
\begin{equation} \labell{table}
 | \l< \xi , \alpha_{j,v} \r> |  <  2 \pi | \l<  y  , \alpha_{j,v} \r> | 
\end{equation}
for each vertex $v$, each edge vector $\alpha_{j,v}$, $j \in I_v$, 
and each $y \in V_\Z^*$ such that $ \l< y , \alpha_{j,v} \r> \neq 0.$
(Each of the sets $\{ \l< y , \alpha_{j,v} \r> \ | \ y \in V_\Z^* \}$
is discrete, so the set of $\xi$'s that satisfy (b) is a neighborhood
of the origin in $V_\C^*$.)
Then
\begin{equation} \labell{GBV}
\sum_{x \in \Delta \cap V_\Z} f(x) = 
\left. \Td_\Delta (\deldel{h_1},\ldots,\deldel{h_d}) \right|_{h=0}
\int_{\Delta(h)} f(x) dx .
\end{equation}
\end{Theorem}

\begin{Remark}
The right hand side of \eqref{GBV} is an infinite sequence.
The theorem asserts that this sequence converges to the left hand side. 
In Appendix \ref{sec:polynomials} we show that this convergence
is uniform on compact subsets of \eqref{table}.
\end{Remark}

The proof of Theorem \ref{EM:exponentials} uses the following 
characterization of $\Gamma_F^\flat$.

\begin{Lemma} \labell{flat characterization}
Let $F$ be a face of $\Delta$.
\begin{enumerate}
\item
If $\gamma \in \Gamma_F$ and $v$ is a vertex 
of $\Delta$ such that $v \in F$, then
$$ e^{2\pi i\l< \gamma,\valpha_{j,v}\r> } = \begin{cases} 
e^{2\pi i\l< \gamma,\valpha_{j,F}\r> } & \text{ for all } j \in I_F \\
 1 & \text{ for all } j \in I_v \ssminus I_F . \end{cases} $$
\item
If $j \in I_F$ and $\gamma \in \Gamma_F^\flat$, then
$e^{2\pi i \l< \gamma , \alpha_{j,F} \r> } \neq 1$.
\item
For $\gamma \in \Gamma_F$, 
$$\gamma \in \Gamma_F^\flat \quad \text{if and only if} \quad
  e^{2 \pi i \l< \gamma, \alpha_{j,F} \r> } \neq 1 
\quad \text{for all} \quad    j \in I_F.$$
\end{enumerate}
\end{Lemma}

\begin{proof}
Let $y \in V_F^* \cap V_\Z^*$ be a representative of $\gamma$
(see \eqref{GammaF}). Then, by definition,
$e^{2\pi i \l< \gamma , \alpha_{j,F} \r> } 
 = e^{2\pi i \l< y , \alpha_{j,F} \r> }$.
Because $y \in V_F^*$, and by \eqref{VFstar}, there exist real numbers
$a_j$, for $j \in I_F$, such that $y = \sum_{j \in I_F} a_j u_j$.
Then
\begin{equation} \labell{for F}
\l< y , \alpha_{j,F} \r> = a_j \qquad \text{ for all } j \in I_F.
\end{equation}
Defining $a_j = 0$ for $j \in I_v \ssminus I_F$, we also have
$y = \sum_{j \in I_v} a_j u_j$, and
\begin{equation} \labell{for v}
\l< y , \alpha_{j,v} \r> = a_j \qquad \text{ for all } j \in I_v.
\end{equation}
In particular,
\begin{equation} \labell{get one}
\l< y , \alpha_{j,v} \r> = 0 \qquad \text{ for all } j \in I_v \ssminus I_F.
\end{equation}
Part (1) follows from \eqref{for F}, \eqref{for v}, and \eqref{get one}.

Fix $j \in I_F$.
Suppose $e^{2\pi i \l< \gamma , \alpha_{j,F} \r> } = 1$.
Then we can choose a representative $y \in V_F^* \cap V_\Z^*$ of $\gamma$
such that $\l< y , \alpha_{j,F} \r> = 0$.
Writing $y = \sum_{l \in I_F} a_l u_l$, we have
$a_j = \l< y , \alpha_{j,F} \r> = 0$.
Let $E$ be the face of $\Delta$ such that $I_E = I_F \ssminus \{ j \}$.
Then $y = \sum_{l \in I_E} a_l u_l$, so, by \eqref{VFstar}, 
$y$ is in $V_E^*$, and so $\gamma \in \Gamma_E$.  
In particular, by \eqref{def:GammaFsharp}, $\gamma \not \in \Gamma_F^\flat$.  
This proves Part (2).

Let $\gamma \in \Gamma_F$.
By Part (2), if $\gamma \in \Gamma_F^\flat$ then 
$e^{2 \pi i \l< \gamma, \alpha_{j,F} \r> } \neq 1$ for all $j \in I_F$.
Conversely, suppose that $\Gamma \not\in \Gamma_F^\flat$.
Then, by \eqref{def:GammaFsharp}, there exists a face $E$ such that
$\gamma \in \Gamma_E$ and $E \supsetneq F$.
Let $j \in I_F \ssminus I_E$.  Let $v$ be any vertex of $F$
(and hence of $E$).  Then
$e^{2\pi i \l< \gamma , \alpha_{j,F} \r> } = 
 e^{2\pi i \l< \gamma , \alpha_{j,v} \r> } = 1$,
where the first equality follows from Part (1) for the face $F$,
and where the second equality follows from Part (1) for the face $E$.
This proves Part (3).
\end{proof}

\begin{Claim*} 
For each $v \in \Vert(\Delta)$,
\begin{equation} \labell{Tdv}
 \Td_\Delta(S_1,\ldots,S_d) = \Td_v(\{ S_j \}_{j \in I_v} )
 + \text{ multiples of } S_j \text{ for } j \not\in I_v .
\end{equation}
\end{Claim*}

\begin{proof}
Recall that
$$ \Td_\Delta(S_1,\ldots,S_d) = \sum_F \sum_{\gamma \in \Gamma_F^\flat}
   e^{ 2 \pi i \left< \gamma , F \right> } 
   \prod_{j \not\in I_F} \frac{S_j}{1-e^{-S_j}}
   \prod_{j \in I_F} \frac{S_j}
        {1 - e^{2 \pi i \left< \gamma,\alpha_{j,F} \right> } e^{-S_j} }.$$
By Part (2) of Lemma \ref{flat characterization},
for each $\gamma \in \Gamma_F^\flat$ and $j \in I_F$,
$$ \frac{S_j}{1 - e^{2 \pi i \left< \gamma,\alpha_{j,F} \right> } 
                  e^{-S_j} }
   = \text{ a multiple of } S_j .$$
Because $v \not\in F$ implies that there exists $j \in I_F$
such that $j \not \in I_v$, and because
$\frac{S_j}{1 - e^{-S_j} } = 1 + \text{ a multiple of } S_j$,
$$ \Td_\Delta(S_1,\ldots,S_d) = 
\sum\limits_{\substack{F \text{ such that } \\ v \in F }}
\sum_{\gamma \in \Gamma_F^\flat} e^{2 \pi i \left< \gamma , F \right> }
 \prod_{j \in I_v \ssminus I_F} \frac{S_j}{1 - e^{-S_j} } 
 \prod_{j \in I_F} 
 \frac{S_j}{1 - e^{2 \pi i \left< \gamma , \alpha_{j,v} \right> } e^{-S_j} }
 + \begin{array}{l} \text{ multiples of } S_j \\ 
   \text{ for } j \not\in I_v.\end{array} $$ 
By \eqref{def:GammaFsharp}, 
$$ \Gamma_v = \bigsqcup\limits_{\substack{F \text{ such that } \\ v \in F }}
   \Gamma_F^\flat .$$
Also, $e^{2 \pi i \left< \gamma ,  F \right> } 
       = e^{2 \pi i \left< \gamma , v \right> }$
whenever $v \in F$.  By this and Part (1) of Lemma \ref{flat characterization},
$$ \Td_\Delta(S_1,\ldots,S_d) = \sum_{\gamma \in \Gamma_v}
   e^{2 \pi i \left< \gamma, v \right> }
   \prod_{j \in I_v} \frac{S_j}
        {1 - e^{2 \pi i \left< \gamma , \alpha_{j,v} \right> } e^{-S_j} }
   + \text{ multiples of } S_j \text{ for } j \not\in I_v. $$
By the definition \eqref{Tdv is} of $\Td_v$, this exatly shows \eqref{Tdv}.
\end{proof}

\smallskip\noindent\textbf{Proof of Theorem \ref{EM:exponentials}
 -- Khovanskii-Pukhlikov approach:}\ 


Let $\Omega \subset \R^d$ be the set of all $h \in \R^d$
that are sufficiently small so that the polytope $\Delta(h)$
has the same combinatorics as the polytope $\Delta$
(i.e., the same subsets $I_F \subset \{1,\ldots,d\}$
correspond to faces).
The vertices of the expanded polytope $\Delta(h)$ are then
\begin{equation} \labell{v of h}
   v(h) = v - \sum_{j \in I_v} h_j \alpha_{j,v} ,
\end{equation}
and we have
\begin{equation} \labell{derivative of vh}
 \deldel{h_j} e^{ \l< \xi , v(h) \r> }
 = - \l< \xi , \alpha_{j,v} \r> e^{\l< \xi , v(h) \r> } 
\end{equation}
for any $j \in I_v$.
Let
$$ \cI_\Delta(h,v) = \int_{\Delta(h)} e^{\l< \xi , x \r> } dx 
\quad \text{ and } \quad
   \cS_\Delta(\xi) = \sum_{x \in \Delta \cap V_\Z} e^{\l< \xi , x \r> }. $$
By \eqref{B1},
\begin{equation} \labell{B1 again}
 \cI_\Delta(h,\xi) = \sum_{v \in \Vert(\Delta)}
                     e^{ \l< \xi , v(h) \r> } \cdot
 \frac{1}{|\Gamma_v|} \prod_{j \in I_v} - \frac{1}{\l< \xi,\alpha_{j,v} \r> }
\end{equation}
for all $h \in \Omega$ and $\xi \in V^*_\C$ that lies outside the 
complex hyperplanes
\begin{equation} \labell{happy}
 \l< \xi , \alpha_{j,v} \r> = 0 \quad , \quad v \in \Vert(\Delta)
 \quad , \quad j \in I_v .
\end{equation}
By \eqref{B2},
\begin{equation} \labell{B2 again}
\cS_\Delta(\xi) = \sum_{v \in \Vert(\Delta)} e^{\l< \xi , v \r> }
   \Td_v ( \left\{ - \l< \xi , \alpha_{j,v} \r> \right\}_{j \in I_v} )
   \cdot \frac{1}{|\Gamma_v|} \prod_{j \in I_v} 
        - \frac{1}{\l< \xi , \alpha_{j,v} \r> } ,  
\end{equation}
for all $\xi \in V^*_\C$ that lie outside the hyperplanes \eqref{funny}.
See Remark \ref{rk:cat}.


By \eqref{Tdv}, \eqref{v of h}, and \eqref{B1 again},
\begin{equation} \labell{star five}
 \Td_\Delta(\deldel{h_1},\ldots,\deldel{h_d}) \cI_\Delta(h,\xi)
 = \sum_{v \in \Vert(\Delta)} \Td_v (\left\{ \deldel{h_j} \right\}_{j \in I_v} )
 e^{\l< \xi,v(h) \r> } \cdot \frac{1}{|\Gamma_v|}
 \prod_{j \in I_v} - \frac{1}{\l< \xi , \alpha_{j,v} \r> }.
\end{equation}
So the partial sums of the series \eqref{star five} are 
\begin{equation} \labell{partial sums}
\sum_v P_{m,v} (\{ - \l< \xi , \alpha_{j,v} \r> \}_{j \in I_v} )
     e^{\l< \xi , v(h) \r> }
    \cdot \frac{1}{|\Gamma_v|} 
          \prod_{j \in I_v} - \frac{1}{\l< \xi , \alpha_{j,v} \r> }
\end{equation}
where $P_{m,v}$ are the Taylor polynomials of $\Td_v$.
By Remark \ref{disk}, $P_{m,v} ( \{ S_j \}_{j \in I_v} )$
converges to $\Td_v ( \{ S_j \}_{j \in I_v} )$
uniformly on compact subsets of the polydisk
\begin{equation} \labell{polydisk} 
\{ S \in \C^{I_v} \ | \ |S_j| < b_{j,v} 
   \text{ for all } j \in I_v \},
\end{equation}
where 
\begin{equation} \labell{bjv}
b_{j,v} = \min\limits_
{\substack{y \in V^*_\Z \\ \l< y,\alpha_{j,v} \r> \neq 0}}
2 \pi \l< y , \alpha_{j,v} \r> . 
\end{equation}
Because the functions
$$ e^{\l< \xi , v\r> } \cdot \frac{1}{|\Gamma_v|} \prod_{j \in I_v}
   - \frac{1}{\l< \xi , \alpha_{j,v} \r> } $$
are (continuous, hence) bounded on compact subsets
of the set of $\xi$'s that lie outside the complex hyperplanes
given by \eqref{happy}, the partial sums \eqref{partial sums}
of the series \eqref{star five} converge to
$$ \sum_{v \in \Vert(\Delta)} 
\Td_v ( \left\{ - \l< \xi , \alpha_{j,v} \r> \right\}_{j \in I_v} )
e^{\l< \xi,v(h) \r> } \cdot \frac{1}{|\Gamma_v|} \prod_{j \in I_v}
-\frac{1}{\l< \xi , \alpha_{j,v} \r> } $$
uniformly on compact subsets of the set 
of $(h,\xi) \in \Omega \times V^*_\C$ such that $\xi$ is outside
the hyperplanes \eqref{happy} and satisfies 
 \begin{equation} \labell{hip hip}
  |\l< \xi , \alpha_{j,v} \r> | < b_{j,v}
 \quad \text{ for all } v \in \Vert(\Delta) \text{ and } j \in I_v.
 \end{equation}

Setting $h=0$, by \eqref{B2 again}, we get
\begin{equation} \labell{ultimate}
 \left. \Td_\Delta(\deldel{h_1}, \ldots, \deldel{h_d}) \right|_{h=0}
 \cI_\Delta(h,\xi) = \cS_\Delta(\xi),
\end{equation}
and that the left hand side of \eqref{ultimate} converges to the
right hand side of \eqref{ultimate} uniformly in $\xi$
on compact subsets of the set of $\xi \in V_\C^*$ that lie
outside the hyperplanes \eqref{happy} and in the set \eqref{hip hip}.

However, the right hand side of \eqref{ultimate} and the partial sums
of the left hand side of \eqref{ultimate} are analytic in $\xi$
for \emph{all} $\xi$.

Recall that,
as a consequence of Cauchy's integral formula, if $g_\nu(\xi)$
is a sequence of complex analytic functions on an open subset $U$
of $\C^n$, $g(\xi)$ is a complex analytic function on $U$,
$g_\nu(\xi)$ converges to $g(\xi)$ in $U \ssminus E$
where $E$ is a locally finite union of complex hyperplanes,
and this convergence is uniform on compact subsets of $U \ssminus E$,
then $g_\nu(\xi)$ converges to $g(\xi)$ for \emph{all} $\xi \in U$,
uniformly on compact subsets of $U$.

It follows that \eqref{ultimate} holds for \emph{all} $\xi \in V^*_\C$
that satisfy \eqref{hip hip}, and, moreover, the left hand side
of \eqref{ultimate} converges to the right hand side uniformly in $\xi$
on compact subsets of \eqref{hip hip}.
This gives the Euler-Maclaurin formula for exponential functions 
$e^{\l< \xi,x \r> }$ for \emph{all} $\xi$ in the set \eqref{hip hip}.
It also shows that the limit on the left hand side of \eqref{ultimate} 
commutes with differentiations with respect to $\xi$.
Applying such differentiations to the left and right sides of \eqref{ultimate},
we get
$$ \left. \Td_{\Delta} (\deldel{h_1},\ldots,\deldel{h_d})\right|_{h=0}
   \int_{\Delta(h)} P(x) e^{\l< \xi , x \r> } dx
  = \sum_{x \in \Delta \cap V_\Z} P(x) e^{\l< \xi,x \r> } dx $$
whenever $P(x)$ is a polynomial and $\xi$ is in the set \eqref{hip hip}.

In particular, for $\xi = 0$, we get the Euler-Maclaurin formula
for polynomials.

\smallskip\noindent\textbf{Proof of Theorem \ref{EM:exponentials}
for polynomial functions -- Brion-Vergne approach:}\ 

The terms in \eqref{B1 again} and \eqref{B2 again} are functions 
of $\xi$ whose products with $\prod_{j,v} \left< \xi , \alpha_{j,v} \right>$
extend to analytic functions near $\xi = 0$.
Comparing the Taylor expansions in $\xi$ of the left and right hand terms
of these products, we get
\begin{equation} \labell{C1}
\int_{\Delta} \frac{ \left< \xi , x \right>^k }{k!} dx
   = \sum_{v \in \Vert(\Delta)} 
   \frac{ \left< \xi , v \right>^{k+n} }{ (k+n)! } 
   \cdot \frac{1}{|\Gamma_v|} \prod_{j \in I_v} 
   - \frac{1}{ \left< \xi , \alpha_{j,v} \right> }
\end{equation}
and
\begin{equation} \labell{C2}
\sum_{x \in\Delta \cap V_\Z} \frac{ \left< \xi , x \right>^k }{k!}
 = \sum_{v \in \Vert(\Delta) } \left( e^{\l< \xi, v \r> }
   \Td_v \left( 
    \left\{ - \left< \xi , \alpha_{j,v} \right> \right\}_{j \in I_v} \right) 
   \right)^{ \l< k+n \r> }
   \cdot \frac{1}{|\Gamma_v|} \prod_{j \in I_v} 
   - \frac{1}{\left< \xi , \alpha_{j,v} \right>} 
\end{equation}
where the superscript $\l< k+n \r> $ denotes the homogeneous term
of degree $k+n$ in $\xi$.

Recall that the vertices of $\Delta(h)$ are
$$ v(h) = v - \sum_{j \in I_v} h_j \alpha_{j,v} .$$
$\left< \xi , v(h) \right>^k$ is a polynomial of degree $k$
in the $h_j$'s that only depends on $h_j$ for $j \in I_v$.
For all $j \in I_v$,
$$ \deldel{h_j} \frac{ \left< \xi , v(h) \right>^k }{k!}
   = - \left< \xi , \alpha_{j,v} \right>
     \frac{ \left< \xi , v(h) \right>^{k-1} }{ (k-1)! } .$$
So for any homogeneous polynomial $T(\cdot)$ of degree $\ell$
in the variables $S_j$, $j \in I_v$,
\begin{equation} \labell{yoo hoo}
 \left. T \left( \{ \deldel{h_j} \}_{j \in I_v} \right) \right|_{h=0}
 \frac{ \left< \xi , v(h) \right>^k }{k!} 
 = T \left( \{ - \left< \xi , \alpha_{j,v} \right> \}_{j \in I_v} \right)
 \frac{ \left< \xi , v(h) \right>^{k-\ell} }{ (k-\ell)! }.
\end{equation}

We have
\begin{multline*}
 \Td_\Delta \left( \deldel{h_1} , \ldots , \deldel{h_d} \right)
   \int_{\Delta(h)} \frac{ \left< \xi , x \right>^k }{k!} dx \\
= \sum_{v \in \Vert(\Delta)}
   \Td_v \left( \left\{ \deldel{h_j} \right\}_{j \in I_v} \right)
   \frac{ \left< \xi , v(h) \right>^{k+n} }{ (k+n)! }
   \cdot \frac{1}{ |\Gamma_v| } \prod_{j \in I_v}
   - \frac{1}{ \left< \xi , \alpha_{j,v} \right> }
   \qquad \text{ by \eqref{C1} and \eqref{Tdv} }  \\
   \qquad = \sum_{v \in \Vert(\Delta)}  \sum_{0 \leq \ell \leq k+n}
   \Td_v^{\left< \ell \right> } 
                \left( \left\{ - \left< \xi , \alpha_{j,v} \right> 
   \right\}_{j \in I_v} \right)
   \frac{ \left< \xi , v(h) \right>^{k+n-\ell} }{ (k+n-\ell)! }
   \cdot \frac{1}{|\Gamma_v|} \prod_{j \in I_v}
   - \frac{1}{ \left< \xi , \alpha_{j,v} \right> }
\qquad \text{ by \eqref{yoo hoo}.} 
\end{multline*}
When $h=0$, by \eqref{C2}, this is equal to
$$ \sum_{x \in \Delta \cap V_\Z} \frac{ \left< \xi , x \right>^k }{k!} .$$


\section{The Stokes formula for polytopes and the Cappell-Shaneson algebra}
\labell{sec:stokes}

Khovanskii and Pukhlikov work with derivatives $\deldel{h_i}$
associated to the ``expansion" $\Delta(h) \subset V$ of the polytope.
In this section we give two results that relate such derivatives
to differential operators on $V$.  The first result, 
Proposition \ref{prop:Stokes}, is the Stokes formula for polynomials. 
The second result, Proposition \ref{prop ddIorbi}, is that integration 
over faces can be replaced by differentiations with respect to 
corresponding $h_j$'s.  
A similar argument appears in \cite[Section 3.6]{BV:partition}.
We use these results to define the ``Cappell-Shaneson algebra",
a formalism used by Cappell and Shaneson to express their formulas.
These results play a key role in relating the Cappell-Shaneson formula
to the original Khovanskii-Pukhlikov formula; we do this in Section
\ref{sec:CS formula}.
The two results can also be used to derive the Euler Maclaurin 
formula for polynomials from the formula for exponentials, as we do 
in Section \ref{sec:polynomials}.

\smallskip\noindent\textbf{The Stokes formula for polytopes.}
\begin{Proposition} \labell{prop:Stokes}
Let $V$ be a vector space with a lattice $V_\Z$.
Normalize Lebesgue measure on $V$ so that the measure of a
fundamental domain for the lattice $V_\Z$ is one.
Let $\Delta$ be a rational polytope in $V$.
Let $u_1,\ldots,u_d$ denote the inward normals to its facets,
normalized so that they are primitive elements of the dual lattice 
$V_\Z^*$.  For any $v \in V$, let $D_v$ denote the directional 
derivative in the direction of $v$.  Then, for any $f \in \Cinf(V)$, 
\begin{equation} \labell{Stokes}
\int_\Delta D_v f = - \sum_{i=1}^d \l< u_i , v \r> \int_{\sigma_i} f .
\end{equation}
\end{Proposition}

\begin{proof}
The formula is an immediate consequence of the general Stokes formula.

Alternatively, it follows directly from
$$ \Delta(h_1+\left< u_1,v\right> , \ldots , h_n+\left< u_n,v\right>)
 = \Delta(h) - v.$$
\end{proof}

\smallskip\noindent\textbf{Integration over faces}

\begin{Proposition} \labell{prop ddIorbi}
Let $\Delta$ be a simple polytope and let $F$ be a face of $\Delta$.
Let $\Delta(h)$ be the expanded polytope and $F(h)$ the corresponding
face of $\Delta(h)$.  (See \eqref{Delta}, \eqref{polytope}, and
\eqref{F(h)}.)
Then, for any smooth function $f \in \Cinf(V)$,
the integral of $f$ on $\Delta(h)$ is a smooth function of $h$
for $h$ near $0$, and
\begin{equation} \labell{ddIorbi}
\int_{F(h)} f
= |\Gamma_F| \prod_{i \in I_F} \deldel{h_i} \int_{\Delta(h)} f
\end{equation}
where $\Gamma_F$ is the finite abelian group associated to the face $F$.
(See \eqref{GammaF}.)
\end{Proposition}

In particular,
\begin{equation} \labell{zero if no face}
\deldel{h_{i_1}} \cdots \deldel{h_{i_k}} \int_{\Delta(h)} f = 0
\quad \text{ if } \sigma_{i_1} \cap \ldots \cap \sigma_{i_k} = \emptyset,
\end{equation}
and
\begin{equation} \labell{ddi}
\deldel{h_i} \int_{\Delta(h)} f = \int_{\sigma_i(h)} f .
\end{equation}

\begin{proof}
Choose a polarizing vector $\xi \in V_\Delta^*$
such that if $v \in \VertDelta$, $v \not\in F$, and $x \in F$, 
then $\left< \xi , x \right>  >  \left< \xi , v \right> $.
(For instance, we may take $\xi'$ such that the restriction of
the linear functional
$\langle \xi' , \cdot  \rangle$ to $\Delta$ 
attains its maximum along the face $F$, and take $\xi$ to be
a perturbation of $\xi'$ which is in $V^*_\Delta$.)
Then the edge vectors that are based at a vertex of $F$
but are not contained in $F$ (that is, $\alpha_{j,v}$ for $v \in F$
and $j \in I_F$) are not flipped in the polarization process
\eqref{polarization}.

Let $P_F$ denote the affine plane generated by the face $F$.
After possibly multiplying $f$ by a cut-off function which is equal 
to one near $F$, we may assume that 
$\langle \xi , x \rangle > \langle \xi , v \rangle $
for every vertex $v$ which is not in $F$
and every $x \in P_F \cap \operatorname{supp}(f)$,
where $\operatorname{supp}(f)$ is the support of $f$.

Then for every vertex $v$ which is not contained in $F$ we have
$$ \int_{P_F \cap \bfC_v^\#} f = 0. $$

Similarly, 
\begin{equation} \labell{no contribution}
 \int_{P_{F(h)} \cap \bfC_v^\#(h)} f = 0 
\qquad \text{if $v \not \in F$ and $h$ is sufficiently small},
\end{equation}
where $P_{F(h)}$ is the affine plane generated by $F(h)$ and
where $\bfC_v^\#(h)$ are the cones that occur in
the polar decomposition
\begin{equation} \labell{polar h}
 {\mathbf 1}_{\Delta(h)}
    = \sum_v (-1)^{|\varphi_v|} {\mathbf 1}_{C_v^\#(h)} .
\end{equation}
From \eqref{no contribution} and \eqref{polar h} we get
$$
 \int_{F(h)} f \; = \; \int_{P_F(h)} f \cdot \bfone_{\Delta(h)}
 \; = \; \int_{P_F(h)} f \cdot \sum_{v \in F} (-1)^{|\varphi_v|} 
   \bfone_{\bfC_v^\#(h)}         
 \; = \; \sum_{v \in F} (-1)^{|\varphi_v|} \int_{P_{F(h)} \cap \bfC_v^\#(h)} f. 
$$

It remains to show that
\begin{equation} \labell{goal remains}
 \int_{P_{F(h)} \cap \bfC_v^\#(h)} f \; = \; 
 \left| \Gamma_F \right| \prod_{i \in I_F} 
 \deldel{h_i} \int_{\bfC_v^\#(h)} f 
\end{equation}
for each $v \in \VertDelta$ such that $v \in F$.
Assume without loss of generality that $h_j = 0$ 
for all $j \in I_v \ssminus I_F$.

Let
$$ T_F = \{ r(x-y) \ | \ x,y \in F \} $$
denote the tangent space to the face $F$.  Consider the affine
change of variable map
$$ \varphi \colon T_F \times \R^{I_F} \to V $$
given by
$$ \varphi(y,t) = v + y + \sum_{j \in I_F} t_j \alpha_{j,v} .$$
Let
$$ F_0^\# = \sum_{j \in I_v \ssminus I_F} \R_+ \alpha_{j,v}^\#
\qquad \text{ and } \qquad 
   \R_+^{I_F} (h) = \prod_{j \in I_F} [-h_j,\infty). $$
The map $\varphi$ sends $F_0^\# \times \R_+^{I_F}$ onto $\ol\bfC_v^\#(h)$
and sends $F_0^\# \times \{ (-h_j)_{j \in I_F} \}$ onto $F_v^\#(h)$.
Lebesgue measure in $ F_0^\# \subset T_F$ is normalized
so that the measure of a fundamental chamber for the lattice 
$T_F \cap V_\Z$ is one.
Clearly,
$$ \int_{F_0^\#} f \left( \varphi(y,-h) \right) dy =
   \prod_{i \in I_F} \deldel{h_i} \int_{F_0^\# \times \R_+^{I_F}(h)}
   f \left( \varphi(y,t) \right) dy dt .$$

To conclude \eqref{goal remains} it remains to show that 
$\left| \det d\varphi \right| = \frac{1}{| \Gamma_F |}$.

The map $d\varphi$ sends the subspace 
$TF \times \{ 0 \} \subset TF \times \R^{I_F}$
to the subspace $TF \subset V$
and respects the lattices in these subspaces.
So its determinant is equal to that of the induced map on quotients.
Recall that $V / TF = V_F$.  The induced map on quotients is the map
$$ \ol{\varphi} \colon \R^{I_F} \to V_F $$
given by $\ol{\varphi} ((t_j)_{j \in I_F}) = 
  \sum_{j \in I_F} t_j \alpha_{j,v}$. Its inverse,
$$ \psi \colon V_F \to \R^{I_F} ,$$
is
$$ \psi(\beta) = ( \l< u_j , \beta \r> )_{j \in I_F}.$$
The dual  $\psi^* \colon \R^{I_F} \to V_F^*$  sends the standard basis
element $e_j$ to $u_j$ for each $j \in I_F$. Finally,
$$
\det d \varphi = \det \ol{\varphi} = (\det \psi)\inv 
 = (\det \psi^*) \inv = [ \span_\Z \{ u_j \} : V_F^* \cap V_\Z^* ]\inv 
 = \left| \Gamma_F \right|\inv, 
$$
as desired.
\end{proof}

By \eqref{ddIorbi}, the formulas of Guillemin \eqref{Guillemin} 
and of Brion-Vergne \eqref{BV} are equivalent.

\smallskip\noindent\textbf{The Cappell-Shaneson algebra}

Let $V$ be a vector space with a lattice $V_\Z$ and $\Delta \subset V$
a simple lattice polytope.
Let $\calD$ denote the ring of infinite order constant coefficient
differential operators on $V$.  Consider the algebra
$\calD[[[\sigma_1],\ldots,[\sigma_d]]]$
of power series in the formal variables $[\sigma_i]$, corresponding 
to the facets, and with coefficients in $\calD$.
A general element of this algebra can be written as
$$ A = \sum_\alpha p_\alpha \prod_{i=1}^d [\sigma_i]^{\alpha(i)} $$
where $p_\alpha \in \calD$ for each
$\alpha \colon \{1,\ldots,d\} \to \Z_{\nonneg}$.
Each such element $A$ defines a linear functional $\int A$
which associates to each polynomial $f$ on $V$ the number
\begin{equation} \labell{def of intA}
\int A (f) :=
\sum_\alpha
\left.
\prod_{i=1}^d \left(\deldel{h_i}\right)^{\alpha(i)}
\right|_{h=\lambda} \int_{\Delta(h)} p_\alpha (f) .
\end{equation}
$\int A (f)$ is also defined for every smooth function $f$
for which the right hand side of \eqref{def of intA} is absolutely convergent. 
By \eqref{zero if no face},
\begin{equation} \labell{vanish1}
\int [\sigma_{i_1}] \cdots [\sigma_{i_k}] = 0
\quad \text{
if \ $\sigma_{i_1} \cap \cdots \cap \sigma_{i_k} = \emptyset$.}
\end{equation}
By \eqref{Stokes} and \eqref{ddi},
\begin{equation} \labell{vanish2}
\int \left(
D_v + \sum_{i=1}^d \l< v , u_i \r> [\sigma_i]
\right) = 0 \quad \text{
for all \ $v \in V$}.
\end{equation}
Following Cappell and Shaneson, we consider the quotient
$Q(\Delta)$ of $\calD[[[\sigma_1],\ldots,[\sigma_d]]]$
by the equivalence relations
\begin{equation} \labell{vanish1 in Q}
 [\sigma_{i_1}] \cdots [\sigma_{i_k}] = 0
\end{equation}
for each set $\sigma_{i_1}, \ldots, \sigma_{i_k}$ of facets
where $\sigma_{i_1} \cap \ldots \cap \sigma_{i_k} = \emptyset$,
and, for each $v \in V$,
\begin{equation} \labell{vanish2 in Q}
 D_v + \sum_{i=1}^d \l< v , u_i \r> [\sigma_i] = 0 .
\end{equation}
We call $Q(\Delta)$ the \emph{Cappell-Shaneson algebra}.
(Cappell and Shaneson denote this algebra $Q(\Sigma)$, where $\Sigma$
is the corresponding fan.)

\begin{Remark}
The quotient of the polynomial algebra $\calD[[\sigma_1],\ldots,[\sigma_d]]$
by the relation \eqref{vanish1 in Q} is the Stanley-Reisner ring 
(the face ring) of $\Delta$ with coefficients in $\calD$.
\end{Remark}

By \eqref{vanish1} and \eqref{vanish2}, the definition of
$\int A (f)$ descends to the Cappell-Shaneson algebra $Q(\Delta)$:
for $T \in Q(\Delta)$ that is represented by
$A \in \calD[[[\sigma_1],\ldots,[\sigma_d]]]$,
we can define
\begin{equation} \labell{def of integral}
\int T(f) = \int A(f).
\end{equation}

\begin{Remark}
Cappell and Shaneson do not consider variations $\Delta(h)$
of the polytope $\Delta$.  They consider
the free $\calD$-module $P(\Delta)$ with basis elements $[F]$
corresponding to the faces $F$ of $\Delta$
and the $\calD$-module map $\rho \colon P(\Delta) \to Q(\Delta)$
defined by
$$ \rho([F]) = | \Gamma_F | \prod_{i \in I_F} [\sigma_i].$$
For $T = \rho(\Omega)$, with $\Omega = \sum_F p_F [F] \in P(\Delta)$,
Equation \eqref{ddIorbi} implies
\begin{equation} \labell{CS def of integral}
\int T(f) = \sum_F \int_F p_F f .
\end{equation}
Cappell and Shaneson \emph{define} $\int T(f)$
by \eqref{CS def of integral}.
Comparing with \ref{def of integral}, we see that this is well defined.
\end{Remark}

Define the \emph{summation functional} by
$$ S(f) = \sum_{\Delta \cap V_\Z} f .$$
Using the Cappell-Shaneson algebra, the Khovanski-Pukhlikov formula 
\eqref{eq:KK} for an integral polytope with non-singular fan reads
$$ S = \int \prod_{i=1}^d \frac{[\sigma_i]}{1 - e^{-[\sigma_i]}} ,$$
and the Guillemin-Brion-Vergne formula for an integral simple polytope
reads
$$
  S = \int \sum_{F} \sum_{\gamma \in \Gamma_F^\flat }
\prod_{j \not \in I_F} \frac{[\sigma_j]}{1 - e^{-[\sigma_j]}}
\prod_{j \in I_F} \frac{[\sigma_j]}
     {1 - e^{2\pi i \l< \gamma, \valpha_{j,F} \r> } e^{-[\sigma_j]}}.
$$

\section{The Cappell-Shaneson formula}
\labell{sec:CS formula}

The main differences between the formulas of Cappell and Shaneson
and those of Khovanskii-Pukhlikov, Guillemin, and Brion-Vergne,
are these:  
\begin{itemize}
\item
The latter authors work with expansions of the polytope.
Cappell and Shaneson work with derivatives of the function on $V$.
\item
Cappell and Shaneson express their formula in terms of what we call
the ``Cappell-Shaneson algebra".
\item
Cappell and Shaneson derive their formula for the sum
$\sum_{\Delta \cap V_\Z} f$ from formulas for the weighted sum
$$ \sum_{\Delta \cap V_\Z}{'} f := 
\sum_F \left( \half  \right)^{\codim F} \sum_{\relint(F) \cap V_\Z} f . $$
\end{itemize}

Cappell and Shaneson's exact formulas for simple lattice polytopes,
when applied to polytopes with non-singular fans, 
become the following formulas:

\begin{Theorem}
Let $\Delta$ be an integral lattice polytope in a vector space $V$
with a lattice $V_\Z$.  Let $f$ be a poylnomial function on $V$.  Then
\begin{equation} \labell{CS1}
 \sum_{\Delta \cap V_\Z} f = \int T (f) ,
\end{equation}
\begin{equation} \labell{CS2}
 \sum_{\interior(\Delta) \cap V_\Z} f = \int \hat{T} (f),
\end{equation}
and
\begin{equation} \labell{CS3}
 \sum_{\Delta \cap V_\Z} f
  - \half \sum_{\del \Delta \cap V_\Z} f
  = \int \half(T+\hat{T}) (f),
\end{equation}
where
$$ T = \sum_{F}
\prod_{i \in I_F} \frac{[\sigma_i]}{2}
\prod_{i \not \in I_F} \frac{[\sigma_i]/2} {\tanh([\sigma_i]/2)} $$
and
$$ \hat{T} = \sum_{F}  (-1)^{\codim F}
\prod_{i \in I_F} \frac{[\sigma_i]}{2}
\prod_{i \not \in I_F} \frac{[\sigma_i]/2} {\tanh([\sigma_i]/2)} .$$
\end{Theorem}

The relation of these formula
to the Khovanskii-Pukhlikov formula goes through 
similar formulas that apply to a polytope with some facets 
removed.  We first need a corresponding polar decomposition.

As before, let $\Delta$ be a polytope with $d$ facets.
For a subset $L \subseteq \{1,\ldots,d\}$, denote by $\Delta^L$
the set obtained by removing from $\Delta$ the facets $\sigma_i$, $i \in L$.
In particular, for $L=\emptyset$, $\Delta^L=\Delta$,
and for $L=\{1,\ldots,d\}$, $\Delta^L=\interior(\Delta)$.

Recall that $\Delta = H_1 \cap \ldots \cap H_d$ where each $H_j$ 
is a half-spaces whose boundary is the affine span of the facet $\sigma_j$.
Let 
$$ H_j^L = \begin{cases}
H_j & \text{ if }j \in L \\
\text{interior}(H_j) & \text{ if } j \not \in L.
\end{cases}$$
Then we have
$$ \Delta^L = \bigcap_j H_j^L .$$
Fix a polarizing vector $\xi$ for $\Delta$. Recall that this determines
a subset $\varphi_v$ of $I_v$ for each vertex $v$.
Consider the cones
$$ C_v^{\#,L} = \bigcap_{j \in I_v} H_{j,v}^{\#,L} $$
where
$$ H_{j,v}^{\#,L} = \begin{cases}
H_j^L 			& \text{ if } j \in I_v \ssminus \varphi_v \\
\left(H_j^L\right)^c 	& \text{ if } j \in \varphi_v .
\end{cases} $$

We have the following polar decomposition for $\Delta^L$:
\begin{Proposition}[Polar decomposition with some facets removed]
\begin{equation} \labell{Lawrence without facets}
 {\mathbf 1}_{\Delta^L} (x)
 = \sum_v (-1)^{|\varphi_v|} {\mathbf 1}_{C_v^{\#,L}} (x) .
\end{equation}
\end{Proposition}

\begin{proof}
We shift the bounding hyperplanes of $H_j$ inward or outward according
to whether $j \in L$ or $j \not \in L$.  That is, we shift by an $h$
which belongs to the set
$$ \Orth^L := \{ (h_1,\ldots,h_d) \mid h_j < 0 \text{ for } j \in L
\text{ and } h_j > 0 \text{ for } j \not \in L. \} $$
We have the pointwise limits
$$ {\mathbf 1}_{\Delta^L}(x)
   = \lim\limits_{\substack{h \to 0 \\ h \in \Orth^L}} 
     {\mathbf 1}_{\Delta(h)}(x) $$
and
$$ {\mathbf 1}_{C_v^{\#,L}}(x)
   = \lim\limits_{\substack{h \to 0 \\ h \in \Orth^L}} 
   {\mathbf 1}_{C_v^{\#}(h)}(x) $$
for all $x$.
The proposition follows immediately from these limits
and from the polar decomposition theorem for $\Delta(h)$.
\end{proof}

We have the following variant of the Khovanskii-Pukhlikov formula
for a polytope with some facets removed:

\begin{Proposition} \labell{KK without facets}
Let $V$ be a vector space with a lattice $V_\Z$.
Let $\Delta$ be a lattice polytope in $V$ with facets 
$\sigma_1,\ldots,\sigma_d$ and $f$ a polynomial function on $V$.
Let $L \subset \{ 1, \ldots, d \}$ be any subset. 
Suppose that $\Delta$ is a polytope with a non-singular fan.  Then
$$ \sum_{\Delta^L \cap V_\Z} f \; = \; \int
\prod_{i \in L} \frac{[\sigma_j] e^{-[\sigma_j]}}{1 - e^{-[\sigma_j]}}
\prod_{i \not\in L} \frac{[\sigma_j]} {1 - e^{-[\sigma_j]}} (f) .$$
\end{Proposition}

\begin{proof}
The proof follows exactly the same lines as the proof of 
the Khovanskii Pukhlikov formula, \eqref{eq:KK},
using the polar decomposition
\eqref{Lawrence without facets} for $\Delta^L$.
We leave the details to the reader.
\end{proof}

We derive the following formula for the weighted sum:

\begin{Proposition} \labell{CP half}
Suppose that $\Delta$ is an integral polytope with a non-singular fan 
and $f$ is a polynomial. Then
\begin{equation} \labell{eq:CP half}
\int \prod_{i=1}^d \frac{[\sigma_i]/2}{\tanh([\sigma_i]/2)} (f)
 \; = \; \sum_{\Delta \cap V_\Z}{'} f.
\end{equation}
\end{Proposition}

\begin{proof}
Consider the left hand side of equation \eqref{eq:CP half}:
\begin{equation} \labell{well}
\int \prod_{j=1}^d \frac{[\sigma_j/2]}{\tanh[\sigma_j/2]}  (f).
\end{equation}
Since
$$ \frac{D/2}{\tanh(D/2)}
 = (D/2) \frac{e^{D/2} + e^{-D/2}}{e^{D/2} - e^{-D/2}}
 = \half (1 + e^{-D}) \frac{D}{1 - e^{-D}},$$
\eqref{well} is equal to
$$\int \prod_{j=1}^d \half \left( 1 + e^{-[\sigma_j]} \right)
\frac{[\sigma_j]}{1 - e^{-[\sigma_j]}} (f).$$
Expanding, this becomes
$$ \left( \half \right)^d \sum_{L \subseteq \{1,\ldots,d\}} \int
\prod_{j \in L} \frac{[\sigma_j] e^{-[\sigma_j]}}{1 - e^{-[\sigma_j]}}
\prod_{j \not\in L} \frac{[\sigma_j]}{1 - e^{-[\sigma_j]}} (f),$$
which, by Proposition \ref{KK without facets}, is equal to
\begin{equation} \labell{well well well}
\left( \half \right)^d \sum_L \sum_{\Delta^L \cap V_\Z} f .
\end{equation}
For each face $F$ of $\Delta$, the relative interior of $F$ is contained
in $\Delta^L$ if and only if $I_F \cap L = \emptyset$.  The number
of subsets $L$ which satisfy this condition is $2^{d-\codim F}$.
Therefore, \eqref{well well well} is equal to
$$ \sum_F \left( \half \right)^{\codim F} \sum_{\interior(F) \cap V_\Z} f,$$
which is the right hand side of \eqref{eq:CP half}.
\end{proof}

Lemma \ref{wted vs nonwted} relates weighted sums
to non-weighted sums.  For a face $F$ of the simple polytope $\Delta$, let
$$ \sum_{\Delta \cap V_\Z}{'} f $$
denote the weighted sum with respect to the affine span of $F$, 
that is,
$$ \sum_{\Delta \cap V_\Z}{'} f = 
   \sum\limits_{\substack{E \\ E \subseteq F}} 
   \left(\half\right)^{\dim F - \dim E} 
   \sum_{\relint(E) \cap V_\Z} f .$$ 
Here, the faces of $F$ are exactly the faces $E$ of $\Delta$
that are contained in $F$, and the exponent $\dim F - \dim E$
is the codimension of $E$ in the affine span of $F$.

\begin{Lemma} \labell{wted vs nonwted}
$$ \sum_{\Delta \cap V_\Z} f = \sum_F \left(\half\right)^{\codim F} 
                               \sum_{F \cap V_\Z}{'} f \ ,$$
and
$$ \sum_{\interior(\Delta) \cap V_\Z} f = \sum_F \left(-\half\right)^{\codim F} 
                               \sum_{F \cap V_\Z}{'} f \ .$$
\end{Lemma}

\begin{proof}

$$ \sum_F \left(\pm\half\right)^{\codim F} \sum_{F \cap V_\Z}{'} f 
  = \sum_F \left(\pm\half\right)^{\codim F}
\sum\limits_{\substack{E \\ E \subseteq F}} \left(\half\right)^{\dim F - \dim E}
\sum_{\relint(E) \cap V_\Z} f$$
\begin{equation} \labell{sub in here}
=  \sum_E \left(\half\right)^{\codim E}
\left( \sum\limits_{\substack{F \\ F \supseteq E}} (\pm 1)^{\codim F} \right)
\sum_{\relint(E) \cap V_\Z} f.
\end{equation}
Because $\Delta$ is simple,
$$ \sum\limits_{\substack{F \\ F \supseteq E}} 1 
   = 2^{\codim E} \quad \text{and} \quad
     \sum\limits_{\substack{F \\ F \supseteq E}} (-1)^{\codim F} 
= \begin{cases} 0 & E \subsetneq \Delta \\
                1 & E = \Delta .  \end{cases}$$
Substituting this in \eqref{sub in here} gives the lemma.
\end{proof}

We are now ready to derive the Cappell-Shaneson formula.

For each face $F$ of $\Delta$,
we have $i \in I_F$ if and only if $\sigma_i \supseteq F$.
For $i \not \in I_F$, the intersection $\sigma_i \cap F$
is either empty or is equal to a face of $\Delta$ which is a facet of $F$.
We denote
$$ \Sigma_F = \{ i \mid \sigma_i \cap F \text{ is a facet of } F \}.$$
Since
$$ \frac{[\sigma_i]/2}{\tanh([\sigma_i]/2)}
 = 1 + \text{ a multiple of } [\sigma_i],$$
and by \eqref{vanish1}, we get
\begin{equation} \labell{T}
 T = \sum_F \prod_{i \in I_F} \frac{[\sigma_i]}{2}
\prod_{i \in \Sigma_F} \frac{[\sigma_i]/2}{\tanh([\sigma_i]/2)} 
\end{equation}
and
\begin{equation} \labell{hat T}
 \hat{T} = \sum_F (-1)^{\codim F}
\prod_{i \in I_F} \frac{[\sigma_i]}{2}
\prod_{i \in \Sigma_F} \frac{[\sigma_i]/2}{\tanh([\sigma_i]/2)} 
\end{equation}
as elements of the Cappell-Shaneson algebra $Q(\Delta)$. 
By \eqref{T} and \eqref{ddIorbi}, 
followed by Proposition \ref{CP half} applied to the face $F$, 
and further followed by Lemma \ref{wted vs nonwted}, we get
$$ \int T(f) = \sum_F \left(\half\right)^{\codim F} \int_F
\prod_{i \in \Sigma_F}
\frac{[\sigma_i]/2}{\tanh([\sigma_i]/2)} (f) 
 = \sum_F \left(\half\right)^{\codim F} \sum_{F \cap V_\Z}{'} f 
 = \sum_{\Delta \cap V_\Z} f.$$
Similarly, from \eqref{hat T} we get
$$ \int \hat{T}(f) = \sum_F \left(-\half\right)^{\codim F} \int_F
\prod_{i \in \Sigma_F}
\frac{[\sigma_i]/2}{\tanh([\sigma_i]/2)} (f) 
 = \sum_F \left(-\half\right)^{\codim F} \sum_{F \cap V_\Z}{'} f
 = \sum_{\interior(\Delta) \cap V_\Z} f .$$
This proves \eqref{CS1} and \eqref{CS2}.
The equality \eqref{CS3} clearly follows from these.

\begin{Remark}
We expect that a similar argument will show that the Cappell-Shaneson
formula for simple polytopes is equivalent to the Guillemin-Brion-Vergne 
formula.
\end{Remark}

\begin{Remark}
Formulas with more general weightings have been developed
in \cite{AW}.
\end{Remark}

\appendix

\section{Relation to remainder formulas}

Another exact formula for polynomial functions on simple polytopes
appeared in our recent paper \cite{KSW2}. 
There we proved an Euler Maclaurin formula \emph{with remainder} 
for simple polytopes and gave estimates on the remainder.  From this 
we deduced an exact formula for polynomials directly, without passing 
through formulas for exponential functions.  Let us describe our
exact formula from \cite{KSW2} in our current notation.

Let $\lambda$ be a complex root of unity, say 
$$ \lambda^N = 1.$$
Define a sequence of functions $Q_{m,\lambda}(x)$ on $\R$ recursively, 
as follows.  For $m=1$, set
$$ Q_{1,\lambda} (x) = \frac{\lambda}{1 - \lambda} \sum_{n \in \Z}
 \lambda^n \bfone_{[n,n+1)} (x) .$$
Given the function $Q_{m-1,\lambda}(x)$, define the function
$Q_{m,\lambda}(x)$ by the conditions
$$ \frac{d}{dx} Q_{m,\lambda} (x) = Q_{m-1,\lambda}(x)
   \quad \text{ and } \quad
   \int_0^N Q_{m,\lambda}(x) dx = 0 .$$
Consider the polynomial
$$ \bfM^{k,\lambda} (S) = 
   \left( \frac12 + \frac{\lambda}{1-\lambda} \right) S
   + Q_{2,\lambda} (0) S^2 + \ldots + Q_{k,\lambda} (0) S^k.$$
Let $V$ be a vector space with a lattice $V_\Z$.
Let 
$$ \Delta = \{ x \ | \ \l< u_i , x \r> + \mu_i \geq 0 
   \ , \ i=1,\ldots,d \} $$
be a simple lattice polytope in $V$, where $u_1,\ldots,u_d \in V^*$
are the normals to the facets of $\Delta$, normalized so that they are 
primitive elements of the lattice $V_\Z^*$. Let
$$ \Delta(h) = \{ x \ | \ \l< u_i , x \r> + \mu_i + h_i \geq 0 
   \ , \ i=1,\ldots,d \} .$$
For a face $F$ of $\Delta$, an element $\gamma$ of $\Gamma_F$,
and an index $1 \leq j \leq d$, let
$$ \lambda_{\gamma,j,F} = \begin{cases}
  e^{ 2\pi i \l< \gamma, \alpha_{j,F} \r> } & j \in I_F \\
  1 & j \not \in I_F ,
\end{cases}$$
and consider the differential operators
$$ \bfM^k_{\gamma,F} = 
   \prod_{j=1}^d \bfM^{k,\lambda_{\gamma,j,F}} (\deldel{h_j}) .$$
Let
$$\sum_{\Delta \cap V_\Z}{'} f := 
  \sum_F (1/2)^{\codim F} \sum_{\relint(F) \cap V_\Z} f ,$$
summing over the faces $F$ of $\Delta$.
Then for any polynomial function $f$ on $\Delta$,
for sufficiently large $k$,
\begin{equation} \labell{exact from remainder}
   \sum_{\Delta \cap V_\Z}{'} f = 
   \sum_F \sum_{\gamma \in \Gamma_F^\flat} \bfM^k_{\gamma,F} 
   \left. \int_{\Delta(h)} f \right|_{h=0} .
\end{equation}

\smallskip

For comparison, the Euler Maclaurin formula \eqref{GBV}
for simple lattice polytopes can be written as
$$ \sum_{\Delta \cap V_\Z} f = \sum_F \sum_{\gamma \in \Gamma_F^\flat}
   \left. \bfT_{\gamma,F} \int_{\Delta(h)} f \right|_{h=0} $$
with
$$ \bfT_{\gamma,F} = \prod_{j=1}^d T^{\lambda_{\gamma,j,F}} 
   \left( \deldel{h_j} \right) 
\qquad \text{ and } \qquad
   \bfT^\lambda(S) = \frac{S}{1 - \lambda e^{-S}} .$$

A similar argument (see below) gives 
\begin{equation} \labell{Thm1prime}
 \sum_{\Delta \cap V_\Z}{'} f = \sum_F \sum_{\gamma \in \Gamma_F^\flat}
 \left. \bfL_{\gamma,F} \int_{\Delta(h)} f \right|_{h=0} 
\end{equation}
with
$$ \bfL_{\gamma,F} = \prod_{j=1}^d \bfL^{\lambda_{\gamma,j,F}} 
   \left( \deldel{h_j} \right) $$
and
$$
 \bfL^\lambda(S) = \frac{S}{2} \cdot 
                     \frac{1 + \lambda e^{-S}}{1 - \lambda e^{-S}} 
 = s \cdot \left( \half + \lambda e^{-S} + \lambda^2 e^{-2S} 
                    + \lambda^3 e^{-3S} + \ldots      \right).$$

As observed by Mich\`{e}le Vergne \cite{V:private}, 
the equivalence of Formulas \eqref{exact from remainder} and \eqref{Thm1prime}
is seen from
\begin{Lemma} \labell{lem:taylor}
$\bfM^{k,\lambda}(S)$ is the $k$th Taylor polynomial of $\bfL^\lambda(S)$.
\end{Lemma}

We complete this section by giving the proofs of \eqref{Thm1prime}
and of Lemma \ref{lem:taylor}.

\begin{proof}[Proof of \eqref{Thm1prime}]
We have the following analogue of \eqref{claim1}:
\begin{equation} \labell{claim1-prime}
\sum_{x \in \bfC_v \cap V_\Z}{'} e^{\l< \xi , x \r> }
\; = \; e^{\l< \xi,v \r> } \cdot \frac{1}{|\Gamma|} \sum_{\gamma \in \Gamma}
e^{2\pi i \l< \gamma,v \r> } \prod_{j=1}^n
\frac12 \cdot 
\frac{1 + e^{2\pi i \l< \gamma,\valpha_j \r> } e^{\l< \xi,\valpha_j \r> } }
     {1 - e^{2\pi i \l< \gamma,\valpha_j \r> } e^{\l< \xi,\valpha_j \r> } }.
\end{equation}
Indeed, applying 
$$ \half \cdot \frac{1 + \lambda e^{-S}}{1 - \lambda e^{-S}}
   = \half + \lambda e^{-S} + \lambda^2 e^{-2S} + \ldots $$
to $\lambda = e^{2 \pi i \l< \gamma,v \r> }$
and $e^{-S} = e^{\l< \xi , \alpha_j \r> }$
and rearranging the terms, the right hand side of \eqref{claim1-prime}
is equal to
$$ \sum_{(k_1,\ldots,k_n) \in \Z_{\geq 0}^n}{'}
   e^{ \l< \xi , v+ \sum k_j \alpha_j \r> }
   \cdot \frac{1}{|\Gamma|} \sum_{\gamma \in \Gamma} 
         e^{2 \pi i \l< \gamma , v + \sum k_j \alpha_j  \r> } $$
which, by \eqref{frobenius}, is equal to the right hand side of
\eqref{claim1-prime}.

From this we get the following analogue of \eqref{BB2}:
\begin{equation} \labell{BB2-prime}
\begin{aligned}
\sum_{x \in {\ol{\bfC}}_v^\sharp \cap V_\Z} {'} e^{\left< \xi , x \right> } 
 &= e^{\left< \xi , v_\shift \right>} \cdot \frac{1}{| \Gamma |} 
    \sum_{\gamma \in \Gamma} e^{ 2 \pi i {\left< \gamma , v \right> } }
    \prod_{j=1}^n  \half
    \frac{1 + e^{ 2 \pi i \left< \gamma , \alpha_j^\sharp \right> }
                 e^{\left< \xi , \alpha_j^\sharp \right> } }
         {1 - e^{ 2 \pi i \left< \gamma , \alpha_j^\sharp \right> }
              e^{\left< \xi , \alpha_j^\sharp \right> } }
 \qquad \text{ by \eqref{claim1-prime} } \\
 &= e^{\left< \xi , v \right> } \cdot \frac{1}{|\Gamma|} 
    \sum_{\gamma \in \Gamma} e^{ 2 \pi i \left< \gamma , v \right> }
    \prod_{j \not\in \varphi_v }  \half
    \frac{1 + e^{2\pi i \left< \gamma , \alpha_j \right> } 
              e^{\left< \xi , \alpha_j \right> } }
         {1 - e^{2\pi i \left< \gamma , \alpha_j \right> } 
              e^{\left< \xi , \alpha_j \right> } }
    \prod_{j \in \varphi_v }  \half
    \frac{1 + e^{-2\pi i \left< \gamma , \alpha_j \right> } 
              e^{-\left< \xi , \alpha_j \right> } }
         {1 - e^{-2\pi i \left< \gamma , \alpha_j \right> } 
              e^{-\left< \xi , \alpha_j \right> } }
 \qquad \text{by \eqref{alpha j sharp}}  \\
 &= (-1)^{|\varphi_v|} e^{\left< \xi , v \right>} \cdot \frac{1}{|\Gamma|}
    \sum_{\gamma \in \Gamma} e^{2 \pi i \left< \gamma , v \right> } 
    \prod \half \frac{1 + e^{2 \pi i \left< \gamma , \alpha_j \right> } 
                          e^{ \left< \xi , \alpha_j \right> } } 
                     {1 - e^{2 \pi i \left< \gamma , \alpha_j \right> } 
                          e^{ \left< \xi , \alpha_j \right> } } \\
 & \qquad \qquad  \text{ by applying the relation } 
   \frac{1+e^x}{1-e^x} = - \frac{1+e^{-x}}{1-e^{-x}} \text{ to } x = 
    - 2 \pi i \left< \gamma,\alpha_j\right> - \left< \xi , \alpha_j \right> 
   \text{ for } j \in \varphi_v.
\end{aligned}
\end{equation}

Let $\bfone^w_\Delta(x)$ denote the weighted characteristic function,
given by $\bfone^w_\Delta(x) = \left( \half \right)^{\codim F}$
if $x$ lies in the relative interior of a face $F$,
and $\bfone^w_\Delta(x) = 0$ if $x \not\in\Delta$.
Define $\bfone^w_C(x)$ in a similar manner whenever $C$ is a convex polyhedral
cone.  We have the following analogue of \eqref{Lawrence}:
\begin{equation} \labell{Lawrence-prime}
\bfone^w_\Delta (x) = \sum_v (-1)^{|\varphi_v|} 
\bfone^w_{\ol{\bfC}_v^\sharp} (x) .
\end{equation}
This can be proved directly (see \cite[section 3]{KSW2}),
or it can deduced from \eqref{Lawrence without facets}
using the formulas
$$ \bfone^w_\Delta(x) = \frac{1}{2^d}\sum_{L \subseteq \{ 1 , \ldots, d \} }
\bfone_{\Delta^L}(x) 
\qquad \text{and} \qquad 
\bfone^w_{\ol{\bfC}_v^\sharp} (x) = \frac{1}{2^d}
   \bfone_{\bfC_v^{\sharp,L}} (x).$$

From this we get the following analogue of \eqref{B2}:
\begin{equation} \labell{B2-prime}
\begin{aligned}
  \sum_{x \in \Delta \cap V_\Z } {'} e^{\left< \xi , x \right> } 
   &= \sum_{v \in \Vert(\Delta)} (-1)^{|\varphi_v|} 
      \sum_{x \in \ol{\bfC}_v^\# \cap V_\Z} {'} e^{ \left< \xi , x \right> }
   & \quad \text{ by \eqref{Lawrence-prime}} \\ 
   &= \sum_{v \in \Vert(\Delta)} 
   e^{\left< \xi , v \right> } \cdot \frac{1}{|\Gamma_v|}
  \sum_{\gamma \in \Gamma_v} e^{2 \pi i \left< \gamma , v \right> }
  \prod_{j \in I_v}  \half \cdot
  \frac{1 + e^{2 \pi i \left< \gamma,\alpha_{j,v} \right> } 
            e^{\left< \xi , \alpha_{j,v} \right> } } 
       {1- e^{2 \pi i \left< \gamma,\alpha_{j,v} \right> } 
            e^{\left< \xi , \alpha_{j,v} \right> } } 
  & \quad \text{ by \eqref{BB2-prime}} \\
 &= \sum_{v \in \Vert(\Delta)} e^{\left< \xi , v \right> }
    \bfL_v \left( \{ - \left< \xi , \alpha_{j,v} \right> \} \right)
    \cdot \frac{1}{|\Gamma_v|} \prod_{j \in I_v} 
                     - \frac{1}{\left< \xi , \alpha_{j,v} \right> }
\end{aligned}
\end{equation}
where \comment{check:}
\begin{equation} \labell{Lv is}
\bfL_v\left( S \right) = 
 \sum_{\gamma \in \Gamma_v} e^{2 \pi i \left< \gamma,v \right>}
 \prod_{ j \in I_v }  \bfL^{\lambda_{\gamma,j,v}}(S_j) .
\end{equation}
Note that $\bfL_v(S)$ is analytic on the polydisk $\{ |S_j| < b_j,
\ j \in I_v \}$ that is described in Remark \eqref{disk}.

The operator that appears in \eqref{Thm1prime} can be written as
$\bfL_\Delta(\deldel{h_1},\ldots,\deldel{h_d})$ where
$$ \bfL_\Delta(S_1,\ldots,S_d) = \sum_F \sum_{\gamma \in \Gamma_F^\flat}
   \prod_{j=1}^d \bfL^{\lambda_{\gamma,j,F}}(S_j).$$
We have the following analogue of \eqref{Tdv}:
for each $v \in \Vert(\Delta)$,
\begin{equation} \labell{Tdv-prime}
 \bfL_\Delta(S_1,\ldots,S_d) = \bfL_v(\{ S_j \}_{j \in I_v} )
 + \text{ multiples of } S_j \text{ for } j \not\in I_v .
\end{equation}
This is shown exactly like \eqref{Tdv}, using the facts that
if $\lambda \neq 1$ then $\bfL^\lambda(S)$ is a multiple of $S$
and if $\lambda = 1$ then $\bfL^\lambda(S) = 1 + $ a multiple of $S$.

By \eqref{v of h} and \eqref{B1 again},
\begin{equation} \labell{Hat}
\int_{\Delta(h)} e^{\l< \xi , x \r> } dx = \sum_{v \in \Vert(\Delta)}
 e^{ \l< \xi , v - \sum_{j \in I_v} h_j \alpha_{j,v}  \r> }
 \cdot \frac{1}{|\Gamma_v|} \prod_{j \in I_v} 
 - \frac{1}{\l< \xi , \alpha_{j,v} \r> } .
\end{equation}
\eqref{Thm1prime} follows from \eqref{Hat} and \eqref{Tdv-prime}
by the same arguments as in the proof of Theorem \ref{EM:exponentials}
in section \ref{sec:EM simple polytope}.
\end{proof}

\begin{proof}[Proof of Lemma \ref{lem:taylor}]
Suppose that $\lambda \neq 1$ and $\lambda^N = 1$.  
An argument similar to those in section \ref{sec:Todd} gives 
\begin{multline} \labell{twisted exact EM for interval}
\left.\left( 
\bfL^{\lambda}\left(\deldel{h_1}\right) 
 + \bfL^{\lambda\inv}\left(\deldel{h_2}\right) 
\right)\right|_{h_1 = h_2 = 0}
\int_{-h_1}^{N+h_2} f(x) dx \\
 = \half f(0) + \lambda f(1) + \lambda^2 f(2) + \ldots
   + \lambda^{N-2} f(N-2) + \lambda^{N-1} f(N-1) + \half f(N)
\end{multline}
for all polynomial functions $f(x)$.  Indeed, direct computation of 
$ \half + \lambda e^\xi + \lambda^2 e^{2\xi}
   + \ldots + \lambda^{N-1} e^{(N-1)\xi} + \half $,
followed by multiplication by $\xi$ and taking the $N$-th degree term 
in the Taylor expansion, gives the following analogue of \eqref{Hoo Hoo}:
\begin{equation} \labell{Hoo Hoo-prime}
 \xi \sum_{x \in [0,N] \cap \Z} {'}(\xi) = 
 \left( \bfL^{\lambda\inv}(\xi) \cdot e^{\xi N} 
      - \bfL^{\lambda}(-\xi) \cdot 1 \right)^{\l< N+1 \r> } 
\end{equation}
whenever $\lambda e^\xi \neq 1$,
where the superscript $\l< N+1 \r> $ denotes the  $N+1$th term
in the Taylor expansion.  
On the other hand, 
\begin{multline*}
 \left. \bfL^{\lambda}(\deldel{h_2}) + \bfL^{\lambda\inv}(\deldel{h_1}) 
 \right|_{h=0}   \xi \int_{-h_2}^{N+h_1} \frac{ (\xi x)^N }{N!} dx  \\
 =
  \left.\bfL^{\lambda\inv}(\deldel{h_1})\right|_{h_1=0}
   \frac{ (\xi(N+h_1))^{N+1} }{ (N+1)! }
   - \left. \bfL^\lambda(\deldel{h_2})\right|_{h_2=0}
   \frac{ (-\xi h_2)^{N+1} }{(N+1)!} \qquad \text{ by \eqref{Hee Hee} } \\
 = \left( \bfL^{\lambda\inv}(\xi) e^{\xi N} - \bfL^\lambda(-\xi) \right)
   ^{ \l< N+1 \r> } \qquad \text{ by \eqref{elephant} } .
\end{multline*}

By direct computation, the first Taylor coefficient of $\bfL^\lambda(S)$
is $\half + \frac{\lambda}{1 - \lambda}$ and that of $\bfL^{\lambda\inv}(S)$
is $\half + \frac{\lambda\inv}{1 - \lambda\inv} = \half - \frac{1}{1-\lambda}$.
Let $a_m$ denote the $m$th Taylor coefficient of $\bfL^\lambda(S)$.
Since $\bfL^\lambda(-S) = \bfL^{\lambda\inv}(S)$, the $m$th Taylor
coefficient of $\bfL^{\lambda\inv}(S)$ is $(-1)^m a_m$.
Taking $F(x)$ to be a polynomial of degree $\leq k+1$ and $f(x) = F'(x)$,
the left hand side of \eqref{twisted exact EM for interval} becomes
\begin{equation} \labell{smily} 
\begin{aligned}
(\half + \frac{\lambda}{1 - \lambda}) f(0)
 + (\half - \frac{1}{1 - \lambda}) f(N)
 + \sum_{m=2}^\infty \left.\left( 
   a_m(\deldel{h_1})^m + (-1)^m a_m (\deldel{h_2})^m 
                     \right)\right|_{h_1 = h_2 = 0}
   \left( F(N) - F(0) \right) \\
 = (\half + \frac{\lambda}{1-\lambda}) f(0)
 + (\half - \frac{1}{1-\lambda}) f(N)
 + \sum_{m=2}^k (-1)^m a_m \left( F^{(m)}(N) - F^{(m)}(0) \right) .
\end{aligned}
\end{equation}
On the other hand, the right hand side of \eqref{twisted exact EM for interval}
is equal to
\begin{multline} \labell{starrr}
   (\half + \frac{\lambda}{1-\lambda}) f(0)
       + (\half - \frac{1}{1-\lambda}) f(N)
   + \frac{\lambda}{1-\lambda} \sum_{n=0}^{N-1} (f(n+1)-f(n)) \\
 = (\half + \frac{\lambda}{1-\lambda}) f(0)
       + (\half - \frac{1}{1-\lambda}) f(N)
   + \frac{\lambda}{1-\lambda} \int_0^N Q_{1,\lambda}(x) f'(x) dx 
   = \ldots \\
= (\half + \frac{\lambda}{1-\lambda}) f(0)
                + (\half - \frac{1}{1-\lambda}) f(N)
  + \sum_{m=2}^k (-1)^m \left. Q_{m,\lambda} (x) f^{(m-1)}(x)\right|_0^N
  + (-1)^{k+1} \int_0^N Q_{k,\lambda} (x) f^{(k)} (x) dx
\end{multline} 
for any $k \geq 2$, by repeated integration by parts
as in the proof of Proposition 27 of \cite{KSW2}.
Recalling that $Q_{m,\lambda}(0) = Q_{m,\lambda}(N)$
and that $\int_0^N Q_{k,\lambda} (x) dx = 0$,
taking $f(x) = F'(x)$ where $F$ is polynomial of degree $\leq k+1$,
the right hand side of \eqref{starrr} becomes 
$$ (\half + \frac{\lambda}{1-\lambda}) f(0)
                + (\half - \frac{1}{1-\lambda}) f(N)
  + \sum_{m=2}^k (-1)^m Q_{m,\lambda} (0) 
    \left( F^{(m)}(N) - F^{(m)}(0) \right) .$$
Comparing this with \eqref{smily} for the monomials
$F(x) = x^3, x^4, x^5, \ldots$ we deduce, by induction on $m$,
that the coefficients $Q_{m,\lambda}(0)$ in $\bfM^{k,\lambda}$  
are equal to the Taylor coefficients $a_m$ of $\bfL^\lambda(S)$
for $m=2,3,4,\ldots$.
\comment{A similar argument should work when $\lambda = 1$.
See our unproved claim in (20) of \cite{KSW2}.}
\end{proof}

\comment{
\begin{enumerate}
\item
We can also apply the twisted Euler Maclaurin formula for a ray 
to the exponential function $f(x) = e^{\xi x}$, $\xi < 0$, 
subtract for the two rays, and use estimates on the remainder.

\item
The fact that $M^{k,\lambda}(S)$ is the $k$th Taylor polynomial of 
$\frac{S}{1 - \lambda^{-S}} - \frac{S}{2}$
was noted in Michele Vergne's email from August.
I (Yael) did not succeed in following her justification of this fact.
She more or less says this:  Let 
$$Q_{0,\lambda}(x) = \sum \lambda^n \delta(x - n).$$
Choose $\Lambda$ such that $e^{2\pi i \Lambda} = \lambda$.
By the Poisson formula, 
$$Q_{0,\lambda} (x) = \sum_{n \in \Z} e^{2 \pi i x (n+\Lambda)}.$$
The $kth$ primitive of this is, up to constants,
$$ Q_{k,\lambda}(x) = \sum_{n \in \Z} 
                   \frac{e^{2 \pi i x (n+\Lambda)}}{(n+\Lambda)^k},$$
At $x=0$,
$$ Q_{k,\lambda}(0) = \sum_{n \in \Z} \frac{1}{(n+\Lambda)^k}.$$
These are related to the Taylor coefficients of
$$ \frac{1}{1 - \lambda \exp(S)} = \frac{1}{1 - \exp(S + 2 \pi i \Lambda)}$$
at $S=0$ because
$$ \frac{1}{\sin^2 x} = \sum_n \frac{1}{(x + n)^2}.$$
\end{enumerate} }

\section{From exponentials to quasi-polynomials: alternative approach}
\labell{sec:polynomials}

One can also deduce the Euler Maclaurin formula
for a polytope directly from an Euler Mclaurin formula for a cone.
This approach is a bit longer than the approaches taken
in Section \ref{sec:EM simple polytope}.  Here we outline this approach
and include a lemma that may be of independent interest.

The exact Euler Maclaurin formula for an exponential function 
on a non-singular convex polyhedral cone is this.
Let $V$ be a vector space with a lattice $V_\Z$,
and let $V_\Z^* \subset V^*$ be the dual lattice.
Let $u_1,\ldots,u_n$ be primitive elements of $V_\Z^*$
which form a basis for $V^*$, and let $\alpha_1,\ldots,\alpha_n \in V$
be the dual basis.  Take any $\lambda_1,\ldots,\lambda_n \in \Z$
and let $v = - \sum_{j=1}^n \lambda_j \alpha_j$.
Consider the finite abelian group
$$ \Gamma = V_\Z^* / \span_\Z \{ u_i \} .$$
Consider the cone
$$ \bfC_v =
\{ x \ | \ \left< u_j , x \right> + \lambda_j \geq 0
 \ , \ j = 1, \ldots , n \}
 = \{ v + \sum t_j \alpha_j \ | \ t_j \geq 0 \ , \ j=1,\ldots,n \} ,$$
and its expansions, given by
$$ \bfC_v(h) =
\{ x \ | \ \left< u_j , x \right> + \lambda_j + h_j \geq 0
 \ , \ j = 1, \ldots , n \} $$
for $h$ near $0$.  Let
$$f(x) = e^{\l< \xi , x \r> },$$
where $\xi \in V^*_\C$ satisfies, for each $j=1,\ldots,n$,
\begin{enumerate}
\item[(a)]
$\Re(\left< \xi , \alpha_j \right> ) < 0$, and
\item[(b)]
$| \l< \xi , \alpha_j \r> | < 2 \pi | \l< y , \alpha_j \r> |$
for all $y \in V_\Z^*$ such that $\l< y , \alpha_j \r> \neq 0$.
\end{enumerate}
(The set of $\xi$'s that satisfy {\rm (b)} is a neighborhood
of the origin in $V^*_\C$.) 
Then
\begin{equation} \labell{calE}
\sum_{\bfC_v \cap V_\Z} f = \left.
 \sum_{\gamma \in \Gamma} e^{2\pi i \l< \gamma,v \r> }
 \prod_{j = 1}^n
 \frac{\deldel{h_j}}
      { 1 - e^{2\pi i \l< \gamma , \valpha_j \r> } e^{-\deldel{h_j}}}
 \right|_{h=0}    \int_{\bfC_v(h)} f .
\end{equation}

This formula follows directly from \eqref{claim1},
\eqref{claim2}, and \eqref{derivative of vh}.

By applying this formula to the polarized cones $\bfC_v^\sharp$
that occur in the polar decomposition (Section \ref{sec:decompose}),
together with some bookkeeping, one deduces the exact Euler Maclaurin
formula on a polytope,
\begin{equation} \labell{Todd Delta}
 \sum_{x \in \Delta \cap V_\Z} e^{\left< \xi, x \right> } =
 \Td_\Delta(\deldel{h_1},\ldots,\deldel{h_d}) \int_{\Delta(h)} 
     e^{\left< \xi , x \right> }  dx
\end{equation}
where $\xi \in V^*$ is sufficiently small and is ``polarizing",
i.e., belongs to the complement of a finite union of (real!) hyperplanes 
through the origin.

One would like to obtain a similar formula for polynomial functions
by taking the derivatives of \eqref{Todd Delta} with respect to $\xi$ 
and taking the limit as $\xi \to 0$,
(or, alternatively, by comparing the coefficients in the Taylor expansions
in $\xi$ of the left and right hand sides of \eqref{Todd Delta}).
For this one needs to show that the infinite order differential operator
$\Td_\Delta(\deldel{h_1}, \ldots, \deldel{h_d})$,
applied to $\int_{\Delta(h)} e^{\left< \xi , x\right> } dx$,
commutes with derivatives and limits with respect to $\xi$.
This follows from the following lemma, which may be of independent interest.

\begin{Lemma} \labell{convergent}
Consider the exponential function
\begin{equation} \labell{f is exp}
 f(\xi,x) = e^{\left< \xi , x \right> } 
\end{equation}
where $x \in V$ and $\xi \in V_\C^*$.
Let $b_1,\ldots,b_d$ be positive numbers, and let $T_\Delta(S_1,\ldots,S_d)$
be a formal power series that converges on the multi-disk 
\begin{equation} \labell{md}
 |S_i| < b_i \quad , \quad  i=1,\ldots,d.
\end{equation}
Then the series 
\begin{equation} \labell{series}
   T_\Delta(\deldel{h_1},\ldots,\deldel{h_d}) 
   \int_{\Delta(h)} f(\xi,x) dx 
\end{equation}
is absolutely convergent whenever $\xi \in V_\C^*$ satisfies the inequalities
\begin{equation} \labell{domain}
 | \left< \xi , \alpha_{i,v} \right> | < b_i  
 \quad \text{ for all $v \in \VertDelta$ and all $i \in I_v$ } ,
\end{equation}
and this convergence is uniform on compact subsets of the domain 
\eqref{domain}.  
\end{Lemma}

\begin{proof}
By general properties of power series, $T_\Delta(S_1,\ldots,S_d)$
is absolutely convergent on the multi-disk \eqref{md}, 
and this convergence is uniform on any strictly smaller multi-disk.  Let
\begin{equation}\labell{Td series}
T_\Delta(S_1,\ldots,S_d) = \sum_F 
\sum\limits_{\substack{\vec{m}=(m_i)_{i\in I_F} \\ 
m_i \text{ non-negative integers} }} 
C_{F,\vec{m}}
    \prod_{i \in I_F} S_i^{1+m_i}  
    + \text{ terms that involve other monomials }.
\end{equation}
The terms that involve other monomials make zero contribution
to \eqref{series}, by \eqref{zero if no face}.

Because the series \eqref{Td series} is absolutely convergent 
on the multi-disk \eqref{md}, so is the series
$$ \sum_F 
\sum\limits_{\substack{ \vec{n} = (n_i)_{i \in I_F} \\ 
n_i \text{ non-negative integers} }}
   \left| C_{F,\vec{n}+\vec{\delta}} \right| 
   \prod_{i \in I_F} S_i^{n_i} $$
for each $\vec{\delta} \in \Z^{I_F}$,
where we set $C_{F,\vec{m}} = 0$ if $m_i < 0$ for some $i \in I_F$.
For $\vec{n} = (n_i)_{i\in I_F}$ and $\vec{m} = (m_i)_{i\in I_F}$ 
we write $\vec{n} \leq \vec{m}$
to mean $n_i \leq m_i$ for all $i \in I_F$, and, for such $\vec{n}$,
we write $| \vec{m} - \vec{n} | = \sum m_i - n_i$.
Then the series
$$ \sum_F \sum_{\vec{n}} \left( 
\sum\limits_{\substack{\vec{m} \text{ such that } \\ 
\vec{n} \leq \vec{m} \text{ and } \\ |\vec{m}-\vec{n}| \leq \dim V }}
|C_{F,\vec{m}}|
\right)
\prod_{i \in I_F} S_i^{n_i}$$
is also absolutely convergent on the multi-disk \eqref{md}.

Let $K$ be any compact subset of the set of $\xi$'s that satisfy 
\eqref{domain}.  Choose positive numbers $\lambda$ and $b_i'$ 
such that $0 < b_i' < b_i$ and $0 < \lambda < 1$ and such that
\begin{equation} \labell{nice}
 | \l< \xi , \alpha_{i,v} \r> | \; < \; \lambda b_i' 
\quad \text{ for all $v \in \VertDelta$ and all $i \in I_v$ } 
\end{equation}
for every $\xi \in K$.
To prove the lemma, we will show that the series \eqref{series}
is dominated on $K$ by a multiple of the converging positive series 
$$ \sum_F \sum_{\vec{n}} \left( 
\sum\limits_{\substack{\vec{m} \text{ such that } \\ 
\vec{n} \leq \vec{m} \text{ and} \\ |\vec{m}-\vec{n}| \leq \dim V }}
|C_{F,\vec{m}}|
\right)
\prod_{i \in I_F} (b_i')^{n_i} .$$

For each face $F$ of $\Delta$ and each $i \in I_F$
we choose an element $\talpha_{i,F} \in V$ to be equal to $\alpha_{i,v}$
for an arbitrary vertex $v \in F$.  Then the elements $\talpha_{i,F}$ 
of $V$ have the following properties.
\begin{enumerate}
\item
\begin{equation} \labell{rep of dual}
 \left< u_l , \talpha_{i,F} \right> = \begin{cases}
 1 & \text{ if } \ l=i \\ 0 & \text{ if } \  l \in I_F \ssminus \{ i \} . 
\end{cases} 
\end{equation}
\item
By \eqref{nice}, for every $\xi \in K$,
$$ | \left< \xi , \talpha_{i,F} \right>|  < \lambda b_i'  \quad 
\text{for each face $F$ and index $i \in I_F$.} $$
\end{enumerate}

The Stokes formula \eqref{Stokes} and part (1) of \eqref{rep of dual} 
imply that for each face $F$ of $\Delta$ and each $i \in I_F$
\begin{equation} \labell{moon}
 \int_{\sigma_i} f \; = \; 
 - \sum_{l \not \in I_F} \left< u_l , \talpha_{i,F} \right> \int_{\sigma_l} f
 - \int_\Delta D_{\talpha_{i,F}} f .
\end{equation}

The exponential function \eqref{f is exp} satisfies
$D_\alpha f = \left< \xi , \alpha \right> f$
for any $\alpha \in V$.
Combining these facts with \eqref{ddi} and \eqref{moon}, we get
\begin{equation} \labell{sun}
 \deldel{h_i} \int_{\Delta(h)} f \; = \; 
 - \sum_{l \not \in I_F} \left< u_l , \talpha_{i,F} \right> 
        \deldel{h_l} \int_{\Delta(h)} f
 - \left< \xi , \talpha_{i,F} \right> \int_{\Delta(h)} f .
\end{equation}
Applying $\prod_{i \in I_F} \deldel{h_i}$ to \eqref{sun},
using the fact that the $\deldel{h_i}$'s commute, 
and applying \eqref{ddIorbi}, we get
$$ \deldel{h_i} \frac{1}{|\Gamma_F|} \int_{F(h)} f \; = \; 
 - \sum\limits_{\substack{l \text{ such that } \\ E:= F \cap \sigma_l \\ 
\text{ satisfies } \emptyset \neq E \subsetneq F}}
\left< u_l , \talpha_{i,F} \right> \frac{1}{|\Gamma_E|} \int_{E(h)} f
- \left< \xi , \talpha_{i,F} \right> \frac{1}{|\Gamma_F|} \int_{F(h)} f .$$

Iterating this formula we get, for any $i_1,\ldots,i_k \in I_F$,
\begin{multline} \labell{iterating}
 \prod_{j=1}^k \deldel{h_{i_j}} \frac{1}{|\Gamma_F|} \int_{F(h)} f \; = \; 
(-1)^k \sum
   \limits_{\substack{l_1,\ldots,l_s \text{ such that } \\
    E_r := F \cap \sigma_{l_1} \cap \ldots \cap \sigma_{l_r} \\
   \text {satisfy} \\ F = E_0 \supsetneq \ldots \supsetneq E_s \neq \emptyset}}
\sum_{1 \leq j_1 < \ldots < j_s \leq k} \\
\left( \prod_{j=1}^{j_1-1} \left< \xi , \talpha_{i_j,F} \right> \right) 
 \left< u_{l_1} , \talpha_{i_{j_1},F} \right> 
\left( \prod_{j=j_1+1}^{j_2-1} \left< \xi , \talpha_{i_j,E_1} \right>  \right)
 \left< u_{l_2} , \talpha_{i_{j_2},E_1} \right>
\cdots \\
\cdots \left< u_{l_s} ,  \talpha_{i_{j_s},E_{s-1}} \right>
\left( \prod_{j=j_s+1}^k \left< \xi , \talpha_{i_j,E_s} \right>  \right)
\frac{1}{|\Gamma_{E_s}|} \int_{E_s(h)} f.
\end{multline}

Let $B \geq 1$ be such that 
$$ | \left< u_l , \talpha_{i,E} \right> | \; \leq \; B
\quad \text{ for all } l, i, \text{ and } E. $$
By this and \eqref{nice},
the term the right hand side of \eqref{iterating}
that corresponds to some fixed $l_1,\ldots,l_s$
and some fixed $j_1,\ldots,j_s$ is bounded by
\begin{equation} \labell{bound1}
 B^s \lambda^{k-s} \prod\limits_{\substack{j=1,\ldots,k \hfill \\ 
 j \not \in \{ j_1,\ldots,j_s \}}} b_{i_j}' 
 \cdot \frac{1}{|\Gamma_{E_s}|} \left| \int_{E_s(h)} f \right| 
\end{equation}
where $E_s = F \cap \sigma_{l_1} \cap \ldots \cap \sigma_{l_s}$.
Let $B_1$ be strictly greater than 
$B^s \lambda^{-s} \frac{1}{|\Gamma_E|} \left| \int_{E}  f \right| $
for all $0 \leq s \leq \dim V$ and all faces $E$ of $\Delta$.  
Then, for $h$ near $0$, the bound \eqref{bound1} is less than 
or equal to
\begin{equation} \labell{bound2}
   B_1 \lambda^k \prod\limits_{\substack{
   j=1,\ldots,k \\ j \not\in \{ j_1,\ldots,j_s \} 
   }} b'_{i_j} .
\end{equation}

Let $m_i$ be the number of times that $i$ occurs in $(i_1, \ldots, i_k)$.
Then the left hand side of \eqref{iterating} can be rewritten as
$$ \prod_{i \in I_F} \left( \deldel{h_i} \right)^{1+m_i} \int_{\Delta(h)} f.$$
Let $n_i$ be the number of times that $i$ occurs among $i_j$ 
for $j \in \{ 1 , \ldots , k \} \ssminus \{ j_1, \ldots ,j_s \}$.
Then the bound \eqref{bound2} can be rewritten as
\begin{equation} \labell{bound3}
B_1 \lambda^k \prod_{i \in I_F} (b_i')^{n_i}.
\end{equation}
Denote $\vec{m} = ( m_i, \ i \in I_F)$ and $\vec{n} = ( n_i, \ i \in I_F)$.
Then $n_i \leq m_i$ for all $i$, which we write as
$\vec{n} \leq \vec{m}$, and $\sum (m_i - n_i) \leq s \leq \dim F$,
which we write as $|\vec{m}-\vec{n}| \leq \dim F$.
Then we can further bound \eqref{bound3} by the following number
which depends on $\vec{m}$ and not on $\vec{n}$:
\begin{equation} \labell{bound4}
B_1 \lambda^k \max\limits_{\substack{\vec{n} \text{ such that } \\
\vec{n} \leq \vec{m} \text{ and } \\ |\vec{m}-\vec{n}|\leq \dim F }}
\prod_{i \in I_F} (b_i')^{n_i}.
\end{equation}

This bound was for the summand of \eqref{iterating}
which corresponds to a fixed choice of $l_1,\ldots,l_s$
and of $j_1,\ldots,j_s$.  The number of possible choices
of $l_1,\ldots,l_s$ is bounded by $d^s$,
which is bounded by $d^{\dim F}$, and, further, by $d^{\dim V}$.
The number of possible choices for $j_1,\ldots,j_s$
is $\binom{k}{s}$, which is bounded by $k^{\dim F}$, and, further,
by $k^{\dim V}$.  It follows that the right hand side of \eqref{iterating}
is bounded by
\begin{equation} \labell{bound5}
d^{\dim V} k^{\dim V} B_1 \lambda^k 
  \max\limits_{\substack{\vec{n} \text{ such that } \\
  \vec{n} \leq \vec{m} \text{ and } \\ |\vec{m}-\vec{n}|\leq \dim F }}
  \prod_{i \in I_F} (b_i')^{n_i}.
\end{equation}

Since 
$\begin{CD} k^{\dim V} \lambda^k 
 @>> k \to \infty > 0
\end{CD}$
and since the maximum among positive numbers is bounded
by their sum,
there exists $B_2$ such that \eqref{bound5} is bounded by
\begin{equation} \labell{bound6}
B_2 \sum\limits_{\substack{\vec{n} \text{ such that } \\
  \vec{n} \leq \vec{m} \text{ and } \\ |\vec{m}-\vec{n}|\leq \dim F }}
  \prod_{i \in I_F} (b_i')^{n_i}.
\end{equation}

To conclude, we found $B_2$ such that
$$ \left| \prod_{i \in I_F} \left( 
   \deldel{h_i} \right)^{1+m_i} \int_{\Delta(h)} f \right| \; \leq \; 
   B_2 \sum\limits_{\substack{\vec{n} \text{ such that } \\
     \vec{n} \leq \vec{m} \text{ and } \\ |\vec{m}-\vec{n}|\leq \dim F }}
     \prod_{i \in I_F} (b_i')^{n_i}. $$

The series \eqref{series} can be written as
$$ \sum_F \sum_{\vec{m}} C_{F,\vec{m}} 
   \prod_{i \in I_F} (\deldel{h_i})^{1+m_i}
   \int_{\Delta(h)} f. $$

By what we have shown, for $\xi$ in the compact set $K$, 
this series is dominated by the positive series
$$ \sum_F \sum_{\vec{m}} |C_{F,\vec{m}}|
   B_2 \sum\limits_{\substack{\vec{n} \text{ such that } \\
   \vec{n} \leq \vec{m} \text{ and } \\ |\vec{m}-\vec{n}|\leq \dim F }}
   \prod_{i \in I_F} (b_i')^{n_i}, $$
which can also be re-written as
$$ B_2 \sum_F \sum_{\vec{n}} \left( 
\sum\limits_{\substack{\vec{m} \text{ such that } \\
\vec{m} \geq \vec{n} \text{ and } \\ |\vec{m}-\vec{n}| \leq \dim V }}
|C_{F,\vec{m}}| \right) \prod_{i \in I_F} (b_i')^{n_i} .$$
As we have shown, this series is convergent.
\end{proof}

\end{document}